\documentclass{amsart}

\usepackage[alphabetic]{amsrefs}
\usepackage{amsmath,amsthm,amssymb}            
\usepackage[all]{xy}   
\usepackage{pstricks,color}

\newtheorem{theorem}{Theorem}[section]

\theoremstyle{definition}
\newtheorem{definition}[theorem]{Definition}
\newtheorem{example}[theorem]{Example}
\newtheorem{problem}[theorem]{Problem}

\theoremstyle{remark}
\newtheorem{remark}[theorem]{Remark}

\numberwithin{equation}{section}

\newcommand{\M}{\mathbb{M}}
\newcommand{\F}{\mathbb{F}}
\newcommand{\C}{\mathbb{C}}
\newcommand{\Q}{\mathbb{Q}}
\newcommand{\R}{\mathbb{R}}
\newcommand{\A}{\mathbb{A}}
\newcommand{\Z}{\mathbb{Z}}

\newcommand{\tmf}{\textit{tmf}}
\newcommand{\mmf}{\textit{mmf}}

\newcommand{\cl}{\mathrm{cl}}
\newcommand{\Top}{\mathrm{top}}

\DeclareMathOperator{\Ext}{Ext}

\DeclareMathOperator{\Sq}{Sq}

\DeclareMathOperator{\Char}{char}

\newcommand{\cirrad}{0.1}
\newgray{gridline}{0.8}

\renewcommand{\smash}{\wedge}

\newcommand{\unit}{\mathbf{1}}
\newcommand{\s}{\mathsf{s}}
\newcommand{\f}{\mathsf{f}}
\newcommand{\eff}{\mathsf{eff}}
\newcommand{\veff}{\mathsf{veff}}
\newcommand{\SH}{\mathbf{SH}}
\newcommand{\PP}{\mathbb{P}}

\newcommand{\EE}{\mathbf{E}}

\newcommand{\K}{\mathbf{K}}
\newcommand{\MGL}{\mathbf{MGL}}
\newcommand{\MBP}{\mathbf{MBP}}
\newcommand{\KGL}{\mathbf{KGL}}
\newcommand{\KQ}{\mathbf{KQ}}
\newcommand{\KW}{\mathbf{KW}}
\newcommand{\kgl}{\mathbf{kgl}}
\newcommand{\kq}{\mathbf{kq}}

\newcommand{\holim}{\text{holim}}

\newcommand{\MU}{\mathbf{MU}}
\newcommand{\BP}{\mathbf{BP}}
\newcommand{\KU}{\mathbf{KU}}
\newcommand{\KO}{\mathbf{KO}}
\newcommand{\ku}{\mathbf{ku}}
\newcommand{\ko}{\mathbf{ko}}

\newcommand{\MM}{\mathbf{M}}

\newcommand{\cotensor}{\mathbin{\Box}}

\begin{document}

\title{Motivic stable homotopy groups}

\author{Daniel C. Isaksen}
\address{Department of Mathematics\\
Wayne State University\\
Detroit, MI 48202, USA}
\email{isaksen@wayne.edu}
\thanks{The author was supported by NSF grant DMS-1202213.}

\author{Paul Arne \O stv\ae r}
\address{Department of Mathematics\\
University of Oslo\\
0316 Oslo, Norway}
\email{paularne@math.uio.no}
\thanks{The author was supported by the RCN Frontier Research Group Project no. 250399 ``Motivic Hopf equations", 
a Friedrich Wilhelm Bessel Research Award from the Alexander von Humboldt Foundation, and a Nelder Visiting Fellowship from Imperial College London.}



\date{}

\begin{abstract}
We survey computations of stable motivic homotopy groups over
various fields.  The main tools are the motivic Adams spectral
sequence, the motivic Adams-Novikov spectral sequence, and the
effective slice spectral sequence.  We state some projects for future
study.
\end{abstract}

\maketitle

\section{Introduction}

Motivic homotopy theory was developed by Morel and Voevodsky 
\cite{Morel99b} \cite{MV99} in the 1990's.  The original motivation for 
the theory was to import homotopical techniques into algebraic geometry.
For example, it allowed for the 
powerful theory of Steenrod operations in algebro-geometric cohomology
theories.
Motivic homotopy theory was essential to the solution of several
long-standing problems in algebraic $K$-theory, such as the
the Milnor conjecture \cite{Voevodsky03b} and the
Bloch-Kato conjecture \cite{Voevodsky11}.

The past twenty
years have witnessed a great expansion of motivic homotopy theory.
Just as in classical homotopy theory, the motivic version of 
stable homotopy groups are among the most fundamental invariants.
The goal of this article is to describe what is known about
motivic stable homotopy groups and to suggest directions
for further study.
This topic is just one of several active directions within motivic
homotopy theory.

There are at least two motivations for studying motivic stable
homotopy groups closely.  First, these groups ought to carry 
interesting arithmetic information about the structure of the
ground field; this is essentially the same as the original
purpose of motivic homotopy theory.  

But another motivation has arisen in recent years.  It turns out that
the richer structure of motivic stable homotopy theory sheds new light on
the structure of the classical stable homotopy groups.  In other words,
motivic stable homotopy theory is of interest in classical homotopy theory,
independently of the applications to algebraic geometry.

The motivic stable homotopy category 
$\SH(k)$ is constructed roughly as follows.
Start with a ground field $k$, and consider the category of smooth
schemes over $k$.  Expand this category to a larger category with better
formal properties.  Finally, impose homotopical relations, especially that
the affine line $\A^1$ is contractible.  The result of this process
is the unstable motivic homotopy category over $k$.

Next, identify a bigraded family $S^{p,q}$ of unstable spheres for
$p \geq q \geq 0$.  These spheres are smash products of the 
topological circle $S^{1,0}$ and the algebraic circle 
$S^{1,1} = \A^1 - 0$.
Then stabilize with respect to this bigraded family to obtain
motivic stable homotopy theory.

A key property of the classical stable homotopy category is that
every spectrum is (up to homotopy) 
built from spheres by homotopy colimits.
The motivic analogue of this property does not hold.  In other words,
there exist motivic spectra that are not 
cellular.
Non-trivial field extensions of the base field are 
examples of non-cellular motivic spectra.
Elliptic curves are another source of such examples.
In both cases, interesting arithmetic properties of the algebraic object
interfere.

One important consequence is that the motivic stable homotopy groups
$\pi_{p,q} X = [S^{p,q}, X]$ do not detect equivalences.  One solution
to this problem is to consider \emph{motivic stable homotopy sheaves}.
The idea is to keep track of maps not only from spheres, but maps out
of all smooth schemes.  By construction of the motivic stable 
homotopy category $\SH(k)$,
these smooth schemes serve as generators for the category.
Therefore, motivic stable homotopy sheaves do detect motivic equivalences.

We will not pursue the motivic stable homotopy sheaf perspective further
in this article.  Instead, we will focus just on the groups
$\pi_{p,q} = [S^{p,q}, S^{0,0}]$.  
The precise relationship is that
the motivic stable homotopy groups are the global
sections of the motivic stable homotopy sheaves.
For cellular motivic spectra, the motivic stable homotopy groups
do detect equivalences, and the most commonly studied motivic
spectra are typically cellular.  
Moreover, the stable homotopy sheaves are ``unramified", which means
that they are actually determined by their
sections over fields.  So a thorough understanding of motivic stable
homotopy groups over arbitrary fields leads back to complete
information about the sheaves as well.

\subsection{Completions}

Classically, Serre's finiteness theorem says that in positive dimensions,
every stable homotopy group is finite.  Therefore, it is enough to
compute the $p$-completions of the stable homotopy groups for all
primes $p$.  Tools such as the Adams spectral sequence and the
Adams-Novikov spectral sequence allow for the calculation of these
$p$-completions.  

As one would expect,
the motivic situation is more involved.
The motivic stable homotopy groups are not finite in general.
For example, the group $\pi_{-1,-1}$ contains a copy of the
multiplicative group $k^\times$ of the base field.  
When $k = \C$, this is an uncountable group that is infinitely divisible.
Therefore, all of the $p$-completions of $\pi_{-1,-1}$ vanish, even though
$\pi_{-1,-1}$ itself is non-zero.
Nevertheless, we will compute $p$-completions, even though
there is a certain loss of information.

A full accounting of the machinery of completion in the motivic
context is beyond the scope of this article.  We refer to
\cite{HKO11} and \cite{DI10}.  The completions arise naturally
when considering convergence for the Adams and Adams-Novikov spectral
sequences.

The situation with completions is in fact even more complicated.  
The first motivic Hopf map $\eta$ is the projection
$\A^2 - 0 \to \mathbb{P}^1$, i.e., a map $S^{3,2} \to S^{2,1}$
representing an element in the motivic stable homotopy group
$\pi_{1,1}$.
We need to complete with respect
to $\eta$ in order to obtain convergence for the relevant spectral sequences.  
The arithmetic square
\[
\xymatrix{
\unit \ar[r] \ar[d] & \unit[\eta^{-1}] \ar[d] \\
\unit^{\wedge}_{\eta} \ar[r] & \unit^{\wedge}_\eta [\eta^{-1}] 
}
\]
is a homotopy pullback for the motivic sphere spectrum
$\unit$, as discussed in 
\cite{roendigs-spitzweck-oestvaer.1line}*{Lemma 3.9}.
Thus, information about the $\eta$-completion 
\index{completion, $\eta$}
\index{$\eta$-completion}
and the 
$\eta$-localization (and how they fit together) leads to
information about the uncompleted motivic stable homotopy groups.

Because of our computational perspective, we will mostly be interested
in torsion in the motivic stable homotopy groups.
The articles \cite{Morel06} and \cite{ALP17} study the rationalized
motivic stable homotopy groups.  After inverting $2$,
the motivic sphere spectrum splits into two summands, and the
rationalizations of these summands can be described in terms of
motivic cohomology with rational coefficients.
We will not pursue this perspective further in this article.

\subsection{Organization}

In Section \ref{sctn:Adams}, we introduce the motivic Adams
spectral sequence, which is one of the key tools for computing
stable motivic homotopy groups.
The discussion includes background on Milnor $K$-theory
of fields, the motivic cohomology of a point, the structure
of the motivic Steenrod algebra, and $\Ext$ groups over
the motivic Steenrod algebra.

In Section \ref{sctn:tmsss},
we introduce the motivic effective slice spectral sequence, 
which
is the other key tool for computing stable motivic homotopy groups.
We explain how the effective slice filtration arises in the stable 
motivic homotopy category $\SH(k)$, 
how to compute the layers of this filtration for the
motivic sphere spectrum, and how to compute the layers of a few
other related motivic spectra.

In Sections \ref{sctn:C} and \ref{sctn:R},
we specialize to the base fields $\C$ and $\R$.  We describe
results derived from the motivic Adams spectral sequence.
This includes calculations of the $\eta$-periodic motivic
stable homotopy groups, significant implications for the
calculation of classical stable homotopy groups,
and exotic nilpotence and periodicity properties.

In Section \ref{sctn:general-field},
we discuss calculations over other fields.
This includes computations over finite fields,
a general discussion of graded commutativity,
Milnor-Witt $K$-theory and its relationship to the
groups $\pi_{n,n}$,
and some results on the groups $\pi_{n+1,n}$.

Section \ref{sctn:other} covers calculations of stable homotopy
groups of other motivic spectra, including the $K$-theory spectra
$\KGL$, $\kgl$, $\KQ$, and $\kq$;
motivic (truncated) Brown-Peterson spectra $\MBP\langle n \rangle$;
and the Witt theory spectrum $\KW$.
The key tool here is the
motivic effective slice spectral sequence. 
We also describe how computations involving $\KW$ are related
to Milnor's conjecture on quadratic forms.
 
Finally, Section \ref{sctn:future} provides
some open-ended projects that may motivate further work
in the subject of stable motivic homotopy groups.

\subsection{Notation}

Our convention for motivic grading is to use the notation $(i,j)$, 
where $i$ is the topological dimension and $j$ is the motivic weight.
Some authors use the notation $(i-j) + j \alpha$, which illuminates
the analogy to equivariant homotopy theory.

The following table of notation is provided for the reader's
convenience. 
Detailed descriptions are provided throughout
the manuscript.

\begin{itemize}
\item
$K^M_*(k)$ is the Milnor $K$-theory of the ground field $k$.
(Section \ref{subsctn:Milnor-K})
\item
$\M_p = H^{*,*}(k; \F_p)$ is the motivic cohomology of the ground
field $k$ with $\F_p$ coefficients. (Theorem \ref{thm:cohlgy-point})
\item
$A_* = A^k_*$ is the dual Steenrod algebra
of motivic cohomology operations. The prime is implicit in this notation.
(Section \ref{subsctn:mot-Steenrod})
\item
$\Ext_k = \Ext_{A}(\M_p, \M_p)$ is the cohomology of the
motivic Steenrod algebra, which serves as the $E_2$-page of the
motivic Adams spectral sequence that converges to the 
completed motivic stable homotopy groups.
(Section \ref{subsctn:A-cohlgy})
\item
$\pi_{i,j} = \pi^k_{i,j}$ is the motivic stable homotopy group over $k$
in degree $(i,j)$, completed at some prime that is implicit
and also at $\eta$.
\item
$\Pi^k_m = \bigoplus_{n \in \Z} \pi^k_{m+n,n}$
is the $m$th Milnor-Witt stem.
Sometimes it is more convenient to study all of the groups in
$\Pi^k_m$ at once, rather than one at a time.
(Sections \ref{subsctn:R-eta-periodic} and \ref{subsctn:1-line})
\item
$K_*^{MW}(k)$ is the Milnor-Witt $K$-theory of a field $k$.
(Section \ref{subsctn:Milnor-Witt})
\item
$\f_q(\EE)$ is the $q$th effective slice cover of a motivic spectrum
$\EE$.
(Section \ref{subsctn:effective-filtration})
\item
$\s_q (\EE)$ is the $q$th layer of a motivic spectrum 
$\EE$ with respect to the effective slice filtration,
i.e., the $q$th slice of $\EE$.
(Section \ref{subsctn:effective-filtration})
\item
$\tilde{\f}_q(\EE)$ is the $q$th very effective slice cover of
a motivic spectrum $\EE$.
(Section \ref{subsctn:effective-filtration})
\item
$\SH(k)$ is the motivic stable homotopy category over the base field $k$.
(Section \ref{sctn:tmsss}).
\item
$\unit = S^{0,0}$ is the motivic sphere spectrum that serves as the
unit object of the stable motivic homotopy category $\SH(k)$.
\item
$S^{i,j}$ is the motivic sphere spectrum of dimension $(i,j)$,
and $\Sigma^{i,j}$ represents the suspension functor given
by smashing with $S^{i,j}$.
\item
$\MM R$ is the motivic Eilenberg-Mac Lane spectrum that represents
motivic cohomology with coefficients in $R$.
Usually $R$ is $\Z$ or $\Z/n$.
(Sections \ref{subsctn:slice-1} and \ref{subsctn:MEMs})
\item
$\KGL$ is the $(2,1)$-periodic algebraic $K$-theory motivic spectrum.
(Sections \ref{subsctn:slice-other} and \ref{subsctn:KGL})
\item
$\kgl = \f_0 \KGL = \tilde{\f}_0 \KGL$ is the $0$th
effective slice cover of $\KGL$, or equivalently the $0$th very effective
slice cover of $\KGL$.
(Sections \ref{subsctn:effective-filtration} and \ref{subsctn:KQ})
\item
$\KQ$ is the $(8,4)$-periodic Hermitian $K$-theory motivic spectrum.
(Sections \ref{subsctn:slice-other} and \ref{subsctn:KQ})
\item
$\kq = \tilde{\f}_0 (\KQ)$ is the $0$th very effective slice cover
of $\KQ$.
(Sections \ref{subsctn:effective-filtration} and \ref{subsctn:KQ})
\item
$\KW$ is the Witt theory spectrum obtained from $\KQ$ by inverting
the first motivic Hopf map $\eta$.
(Sections \ref{subsctn:slice-other} and \ref{subsctn:HWT})
\item
$\MGL$ is the motivic algebraic cobordism spectrum, and $\MU$
is the analogous classical complex cobordism spectrum.
(Section \ref{subsctn:slice-1})
\item
$\MBP$ is the motivic version of the Brown-Peterson spectrum,
and $\MBP\langle n \rangle$ is the truncated version, while
$\BP$ and $\BP \langle n \rangle$ are the classical
analogues.
These constructions depend on the choice of some prime,
which is implicit in the notation.
(Sections \ref{subsctn:slice-1} and \ref{subsctn:MBP})
\end{itemize}

\subsection{Acknolwedgements}

The authors appreciate constructive feedback from
J.\ Hornbostel, 
O.\ R{\"o}ndigs, M.\ Spitzweck, and G.\ Wilson.
\index{R{\"o}ndigs, O.}
\index{Spitzweck, M.}
\index{Wilson, G.}
\index{Hornbostel, J.}

\section{The motivic Adams spectral sequence}
\label{sctn:Adams}

Our techniques for computing motivic stable homotopy groups are
fundamentally cohomological in nature.  So our discussion begins
with the computational properties of motivic cohomology.

\subsection{Milnor $K$-theory and the cohomology of a point}
\label{subsctn:Milnor-K}

Let $k$ be a field.
The \emph{Milnor $K$-theory} 
\index{Milnor $K$-theory}
\index{$K$-theory, Milnor}
$K^M_*(k)$ of $k$ is a graded ring
constructed as follows \cite{Milnor69}.  Let $K^M_0(k)$ be $\Z$, and let
$K^M_1(k)$ be the abelian group $k^\times$.
This presents a notational confusion, because the group structure
on $k^\times$ is multiplicative, but we would like to think
of $k^\times$ as an additive group.  We use the symbol
$[a]$ to represent the element of $K^M_1(k)$ corresponding to the
element $a$ of $k^\times$.  Then we have identities such as
\[
[a b] = [a] + [b].
\]

Now $K^M_*(k)$ is defined to be the graded commutative ring
generated by the elements of $K^M_1(k)$ subject to the 
\emph{Steinberg relations}
\[
[a] \cdot [1-a] = 0
\]
\index{Steinberg relation}
in $K^M_2(k)$ for all $a$ in $k - \{ 0, 1 \}$.

The Milnor $K$-theory of familiar fields can be wildly complicated.
For example, $K^M_1(\C)$ is an uncountable infinitely divisible abelian
group.  Fortunately, we will usually consider the much simpler
Milnor $K$-theory modulo $p$, i.e., $K^M_*(k) / p$, for a fixed prime
$p$.

We give several standard examples of Milnor $K$-theory computations.
More details can be found in \cite{Milnor69}, \cite{Weibel13}*{Section III.7}, or \cite{Magurn02}*{Chapter 14E}.

\begin{example}
$K^M_*(\C) / p = \F_p$, concentrated in degree $0$,
since every element of $\C^\times$ is a $p$th power of some other
element.
\end{example}

\begin{example}
$K^M_*(\R) / 2 = \F_2 [\rho]$,
where $\rho$ is another name for the element $[-1]$ in $K^M_1(\R)$.
This calculation follows from the observation that every non-negative
real number is a square.
On the other hand,
$K^M_*(\R) / p = \F_p$ if $p$ is odd, since
every element of $\R^\times$ is a $p$th power of some other element.
\end{example}

In general, let $\rho$ be the element $[-1]$ in $K^M_1(k)/2$
for any field $k$.  Note that $\rho$ is zero if and only if
$-1$ is a square in $k$.
We will see later in Section \ref{subsctn:mot-Steenrod}
that $\rho$ plays a central role in the structure of the
motivic Steenrod algebra.

\begin{example}
If $p \equiv 1 \bmod 4$, then
$K^M_*(\Q_p) / 2 = \F_2[\pi, u] / (\pi^2, u^2)$,
where $\pi = [p]$ and $u = [a]$ for an element $a$ in $\Q_p^\times$
that maps to a non-square in $\F_p^\times$.
\end{example}

\begin{example}
If $p \equiv 3 \bmod 4$, then
$K^M_*(\Q_p) / 2 = \F_2[\pi, \rho] / (\rho^2, \rho \pi + \pi^2)$,
where $\pi = [p]$ and $\rho = [-1]$.
\end{example}

\begin{example}
$K^M_*(\Q_2) / 2 = \F_2 [ \pi, \rho, u ] / (\rho^3, u^2, \pi^2,
\rho u, \rho \pi, \rho^2 + u \pi)$,
where 
$\pi = [2]$, $\rho = [-1]$, and $u = [5]$.
\end{example}

\begin{example}
\label{ex:KM-Fq}
$K^M_*(\F_q)$ equals $\Z[u]/u^2$, where $u = [a]$
for any generator $a$ of $\F_q^\times$.
\end{example}

For us, the point of these Milnor $K$-theory calculations is that
we can use them describe the motivic cohomology of a point.
Motivic cohomology is bigraded, where the first grading corresponds
to the classical topological degree, and the second corresponds
to the motivic weight.
If $k$ contains a primitive $p$th root of unity $\zeta_{p}$,
we let $\tau$ denote the corresponding
generator of $H^{0,1}(k;\F_p)\cong \Z/p$.
\index{motivic cohomology of a point}

\begin{theorem}
\label{thm:cohlgy-point}
\cite{Voevodsky03b} \cite{Voevodsky11}
Suppose that $p$ and $\Char(k)$ are coprime, and 
that $k$ contains a primitive $p$th root of unity.
Then
the motivic cohomology
$\M_p = H^{*,*}(k; \F_p)$ with coefficients in $\F_p$ is isomorphic to
\[
\frac{K^M_*(k)}{p} [ \tau ],
\]
where $K^M_n(k)/p$ has degree $(n,n)$, and
$\tau$ has degree $(0,1)$.
\end{theorem}

\begin{figure}

\psset{unit=1cm}
\begin{pspicture}(-5,-5)(2,2)

\psgrid[unit=1,gridcolor=gridline,subgriddiv=0,gridlabelcolor=white](-4,-4)(1,1)

\pscircle*(0,0){\cirrad}
\pscircle*(-1,-1){\cirrad}
\pscircle*(-1,-1){\cirrad}
\pscircle*(-2,-2){\cirrad}
\pscircle*(-3,-3){\cirrad}
\pscircle*(-4,-4){\cirrad}

\pscircle*(0,-1){\cirrad}
\pscircle*(0,-2){\cirrad}
\pscircle*(0,-3){\cirrad}
\pscircle*(0,-4){\cirrad}

\pscircle*(-1,-2){\cirrad}
\pscircle*(-1,-3){\cirrad}
\pscircle*(-1,-4){\cirrad}

\pscircle*(-2,-3){\cirrad}
\pscircle*(-2,-4){\cirrad}

\pscircle*(-3,-4){\cirrad}

\psline{->}(0,0)(-4.5,-4.5)
\psline{->}(0,-1)(-3.5,-4.5)
\psline{->}(0,-2)(-2.5,-4.5)
\psline{->}(0,-3)(-1.5,-4.5)
\psline{->}(0,-4)(-0.5,-4.5)

\psline{->}(0,0)(0,-4.5)
\psline{->}(-1,-1)(-1,-4.5)
\psline{->}(-2,-2)(-2,-4.5)
\psline{->}(-3,-3)(-3,-4.5)
\psline{->}(-4,-4)(-4,-4.5)

\rput(1,-5){1}
\rput(0,-5){0}
\rput(-1,-5){-1}
\rput(-2,-5){-2}
\rput(-3,-5){-3}
\rput(-4,-5){-4}

\rput(-5,1){1}
\rput(-5,0){0}
\rput(-5,-1){-1}
\rput(-5,-2){-2}
\rput(-5,-3){-3}
\rput(-5,-4){-4}

\uput{0.2}[120](0,0){$K^M_0/p$}
\uput{0.2}[120](-1,-1){$K^M_1/p$}
\uput{0.2}[120](-2,-2){$K^M_2/p$}
\uput{0.2}[120](-3,-3){$K^M_3/p$}
\uput{0.2}[120](-4,-4){$K^M_4/p$}

\uput{0.2}[30](0,-1){$\tau$}

\end{pspicture}

\caption{
\label{fig:cohlgy-pt}
The cohomology of a point.  The grading is homological to 
align with later homotopy group calculations.  The topological degree
is on the $x$-axis, and the motivic weight is on the $y$-axis.
Each dot in column $-n$ represents a copy of 
$K_n^M(k)/p$.
Vertical lines indicate multiplication by $\tau$.
Diagonal lines connect the source and target of multiplications by elements of $K_1^M(k)/p$.
}
\end{figure}
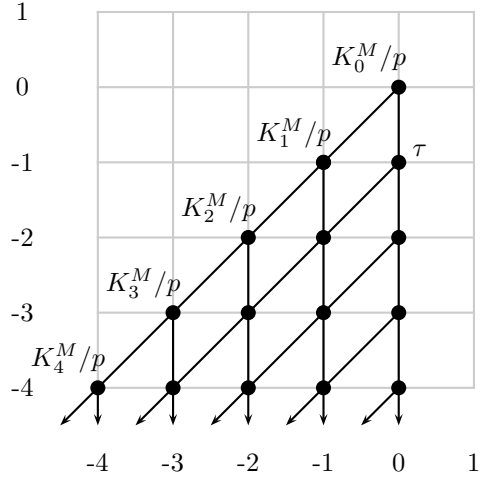

Figure \ref{fig:cohlgy-pt} graphically depicts
the calculation in Theorem \ref{thm:cohlgy-point}.

\subsection{The motivic Steenrod algebra}
\label{subsctn:mot-Steenrod}

Our next task is to record the structure of the motivic Steenrod algebra \cite{MR3730515}, \cite{Voevodsky03a}.
We continue to assume that the base field $k$ contains
a primitive $p$th root of unity and that $\Char(k)$ and $p$ 
are coprime, so that we rely on the calculation
of Theorem \ref{thm:cohlgy-point}.

In classical topology, the dual Steenrod algebra is easier to describe
than the Steenrod algebra itself.  This asymmetry arises from the fact
that the coproduct structure
(i.e., the Cartan formula) is simpler than the 
product structure (i.e., the Adem relations).  

In the motivic context, the advantages of the dual Steenrod algebra
are even more pronounced.  The motivic Adem relations possess
some non-trivial complications that make them difficult to even
write down correctly.  
Beware that more than one published version of the motivic
Adem relations has mistakes.  
See \cite{MR3730515}*{Theorem 5.1} for the correctly formulated
relations.

Let $p$ be an odd prime.
The dual motivic Steenrod algebra $A_*$ is 
\index{dual motivic Steenrod algebra}
\index{motivic Steenrod algebra, dual}
\index{Steenrod algebra, dual motivic}
a commutative
$\M_p$-algebra,
where $\M_p = H^{*,*}(k; \F_p)$ is the cohomology of a point
(see Theorem \ref{thm:cohlgy-point}).
It is 
generated by elements $\tau_i$ for $i \geq 0$ 
and $\xi_j$ for $j \geq 1$, subject to the relations
$\tau_i^2 = 0$.
Here $\tau_i$ has degree $(2 p^i - 1, p^i - 1)$, and
$\xi_j$ has degree $( 2 p^j - 2 , p^j - 1 )$.
The coproduct $\Delta$ on $A_*$ is described by
\[
\Delta (\tau_i) = \tau_i \otimes 1 + 
\sum_{k=0}^i \xi_{i-k}^{p^k} \otimes \tau_k
\]
and
\[
\Delta(\xi_j) = \sum_{k=0}^j \xi_{j-k}^{p^k} \otimes \xi_k.
\]
In order to properly interpret these formulas, we adopt the usual
convention that $\xi_0 = 1$.

The structure of the motivic Steenrod algebra at odd primes
is essentially the same as the structure of the classical Steenrod
algebra.  In fact, in the odd primary case, the motivic Steenrod
algebra is the classical Steenrod algebra tensored over $\F_p$
with a larger coefficient ring.
Consequently, motivic $\Ext$ calculations
over $A_*$ are essentially identical to classical $\Ext$ calculations.

The situation at $p = 2$ is more interesting.
In this case, the dual motivic Steenrod algebra $A_*$ is the commutative
$\M_2$-algebra generated by elements $\tau_i$ for $i \geq 0$ 
and $\xi_j$ for $j \geq 1$, subject to the relations
\[
\tau_i^2 = \tau \xi_{i+1} + \rho \tau_{i+1} + \rho \tau_0 \xi_{i+1}.
\]
Here $\tau_i$ has degree $(2^{i+1} - 1, 2^i - 1)$, and
$\xi_j$ has degree $( 2^{j+1} - 2, 2^j - 1)$.
The coproduct $\Delta$ on $A_*$ is described by
\[
\Delta (\tau_i) = \tau_i \otimes 1 + 
\sum_{k=0}^i \xi_{i-k}^{2^k} \otimes \tau_k
\]
and
\[
\Delta(\xi_j) = \sum_{k=0}^j \xi_{j-k}^{2^k} \otimes \xi_k.
\]
As before, we adopt the usual convention that $\xi_0 = 1$.

\begin{remark}
\label{rem:A-invert-tau}
If we invert $\tau$ and set $\rho$ equal to zero, then we obtain
the classical Steenrod algebra tensored over $\F_2$ with
$\M_2 [\tau^{-1}]$.
The point is that $\xi_{i+1}$ is no longer needed as a generator
because $\xi_{i+1} = \tau^{-1} \tau_i^2$.
Then the elements $\tau_i$ are generators, with no relations.
\end{remark}

Unlike in the classical case,
the motivic Steenrod algebra is more accurately a Hopf algebroid
(in the sense of \cite{Ravenel86}*{Appendix 1}, for example).
For sake of tradition, we will not use the term ``algebroid" in this
context.
This means that the motivic Steenrod
algebra acts non-trivially on the motivic cohomology of a point
in general.
In other words, the right unit $\eta_R: \M_p \to A_*$ is not the
same as the left unit $\eta_L: \M_p \to A_*$.  More precisely,
$\eta_L$ is the obvious inclusion, while
$\eta_R(\tau) = \tau + \rho \tau_0$ and
$\eta_R(\rho) = \rho$.
This complication occurs only if $\rho = [-1]$ is
non-zero in $K^M_*(k)/p$, i.e., only if $-1$ is not a $p$th power
in the ground field $k$.

For definiteness, we restate the action of the Steenrod 
squares on $\M_p$ in more concrete terms.  
At $p=2$, we have $\Sq^1(\tau) = \rho$.  There are 
also non-zero squaring operations on other elements, but they can all
be derived from this basic formula.  

The motivic Cartan formula takes the form
\index{Cartan formula, motivic}
\index{motivic Cartan formula}
\[
\Sq^{2k}(xy) = \sum_{a+b = 2k} \tau^{\epsilon} \Sq^a(x) \Sq^b(y),
\]
where $\epsilon$ is $0$ if $a$ and $b$ are even, while
$\epsilon$ is $1$ if $a$ and $b$ are odd 
\cite{Voevodsky03a}*{Proposition 9.7}.
(There is also a motivic Cartan formula for odd squares, but it
has some additional complications involving $\rho$.)
This implies that
\[
\Sq^2(\tau^2) = \tau \Sq^1(\tau) \Sq^1(\tau) = \rho^2 \tau.
\]

At odd primes, we have that $\beta(\tau) = \rho$, where $\beta$
is the Bockstein operation.  

\subsection{The cohomology of the Steenrod algebra}
\label{subsctn:A-cohlgy}

The next step is to compute the cohomology of the
motivic Steenrod algebra, i.e.,
\[
\Ext_k = \Ext_{A}(\M_p, \M_p).
\]
\index{cohomology of the motivic Steenrod algebra}
\index{motivic Steenrod algebra, cohomology of}
\index{Steenrod algebra, motivic, cohomology of}
This object is a trigraded commutative ring with higher structure
in the form of Massey products.  Two of the gradings correspond to the
degrees familiar in the classical situation.  The additional grading
corresponds to the motivic weight, which has no classical analogue.

The cobar resolution is the standard
way to compute $\Ext$ groups.
\index{cobar resolution}
This works just as well motivically as it does classically.
The cobar resolution is useful for (very) low-dimensional
explicit calculations and for general properties.
In a larger range, one must attack the calculation with
more sophisticated tools such as the May spectral sequence
\cite{May65} \cite{DI10}.

In the modern era, the best way to determine $\Ext$ in a 
range is by computer.
These entirely algebraic invariants are effectively computable in
a large range.  
Computer algorithms typically rely on minimal resolutions rather
than the cobar resolution because it grows more slowly.
In practice, computer calculations of $\Ext$
far outpace our ability to interpret the data with
the Adams spectral sequence 
\cite{Bruner89}
\cite{Bruner91} \cite{Bruner97} \cite{Nassau}
\cite{IWX18}.

Of course, explicit calculations over a field $k$ depend on
specific knowledge of $K^M_*(k)$.

\begin{remark}
\label{rem:A-cohlgy-C}
Over $\C$ and at $p=2$, the cohomology of the motivic Steenrod algebra was
computed with the motivic May spectral sequence through 
the 70-stem \cite{Isaksen14c}.
\index{motivic May spectral sequence}
\index{May spectral sequence, motivic}
\index{spectral sequence, motivic May}
The most recent computer calculations over $\C$ extend far beyond the 100-stem,
with ongoing progress into even higher stems \cite{IWX18}.
\end{remark}

\begin{remark}
Over $\F_q$ and at $p=2$, 
the cohomology of the motivic Steenrod algebra is studied
in detail through the 20-stem \cite{WO17}.  The Milnor $K$-theory
$K^M_*(\F_q)$ is relatively simple, as described in Example \ref{ex:KM-Fq};
this makes explicit computations practical.

If $\Char \F_q \equiv 1 \bmod 4$, then
$-1$ is a square and $\rho$ is zero in $\M_2$.
Therefore, the situation is essentially identical to the calculation
over $\C$.

However, if $\Char \F_q \equiv 3 \bmod 4$,
then $\rho$ is non-zero in $\M_2$.
This case is more complicated (and more interesting).
\end{remark}

\subsection{The $\rho$-Bockstein spectral sequence}
\label{subsctn:rho-Bockstein}

Over $\R$, the cohomology $\M_2$ of a point is $\F_2[\tau, \rho]$.
As discussed in Section \ref{subsctn:mot-Steenrod},
the motivic Steenrod algebra acts non-trivially on a point.
This complication significantly increases the difficulty of
$\Ext$ calculations.
In other terms,
the relation
$\tau_i^2 = \tau \xi_{i+1} + \rho \tau_{i+1} + \rho \tau_0 \xi_{i+1}$
is the source of these difficulties.

In this type of situation, the standard approach is to impose a filtration
that hides these complications (or rather, pushes them into higher filtration).
Then one obtains a spectral sequence that computes $\Ext$.

In this case, we filter by powers of $\rho$ and obtain a
$\rho$-Bockstein spectral sequence that converges to $\Ext_\R$
\cite{Hill11} \cite{DI17}.
\index{$\rho$-Bockstein spectral sequence}
\index{Bockstein spectral sequence, $\rho$}
\index{spectral sequence, $\rho$-Bockstein}
The associated graded object of $A^\R_*$ is easily identified with
$\A^\C_* [\rho]$.
From this observation, one can deduce that
the spectral sequence takes the form
\[
\Ext_\C [ \rho ] \Rightarrow \Ext_\R.
\]

In order for the $\rho$-Bockstein spectral sequence to be practical,
one needs a method for computing differentials.  
The first step in this program is to consider the effect of
inverting $\rho$ \cite{DI17}*{Theorem 4.1}.
If one inverts $\rho$ in $A_*$, then
the relation
\[
\tau_i^2 = \tau \xi_{i+1} + \rho \tau_{i+1} + \rho \tau_0 \xi_{i+1}
\]
can be rewritten in the form
\[
\tau_{i+1} = \rho^{-1} \tau_i^2 + \rho^{-1} \tau \xi_{i+1} + \tau_0 \xi_{i+1}.
\]
This shows that $A_*[\rho^{-1}]$ can be simply described as
$\M_2[\rho^{-1}][\tau_0, \xi_1, \xi_2, \ldots ]$.
The cohomology of $A_*[\rho^{-1}]$ is then straightforward to compute.
One obtains the isomorphism
\[
\Ext_\R [\rho^{-1}] \cong \Ext_{\cl} [\rho^{\pm 1}].
\]
In other words, the $\rho$-periodic part of $\Ext_\R$
is identified with classical $\Ext$ groups.
Beware that this isomorphism 
changes degrees.  For example,
the motivic element $h_{i+1}$ in the $(2^{i+1}-1)$-stem
maps to the classical element $h_i$ in the $(2^i-1)$-stem.

\begin{remark}
In fact, this computational observation extends to an equivalence
of homotopy categories.  Namely, the $\rho$-periodic
$\R$-motivic stable homotopy category 
$\SH(\R)[\rho^{-1}]$
is equivalent to the classical
stable homotopy category 
\cite{Bachmann18}.

Inverting $\rho$ on an $\R$-motivic spectrum is analogous to
taking the geometric fixed points of a $C_2$-equivariant spectrum.
\end{remark}

It remains to compute the $\rho^r$-torsion in $\Ext_\R$ for all $r$,
i.e., to compute the $\rho$-Bockstein $d_r$ differentials for all $r$.
It turns out that these differentials are forced by combinatorial
considerations in a large range.
Every element $x$ of the $\rho$-Bockstein $E_1$-page falls
into one of three categories:
\begin{enumerate}
\item
$x$ survives, 
and $\rho^k x$ is not hit by any differential for all $k$.
These elements are completely known because of the 
calculation of $\Ext_\R [\rho^{-1}]$.
\item
$x$ survives and $\rho^k x$ is hit by some differential for some
$k \geq 1$.
\item
$x$ does not survive.
\end{enumerate}
Somewhat surprisingly, these considerations
determine many $\rho$-Bockstein differentials.
In other words, in a large range, 
there is only one pattern of differentials
that is consistent with sorting all elements into
these three classes.

See \cite{DI17}*{Section 5} for a description of this process in
low dimensions.  Preliminary results indicate that the same naive
method works into much higher dimensions.

\begin{remark}
The $\rho$-Bockstein spectral sequence is useful for fields other 
than $\R$, in which $-1$ is not a square and therefore $\rho$ is non-zero.
Some examples are finite fields $\F_q$ for which $\Char \F_q \equiv 3 \bmod 4$ \cite{WO17}.
\end{remark}

\subsection{The motivic Adams spectral sequence}
\label{subsctn:mot-Adams}

The point of computing $\Ext_k$ is that it serves
as the input to the motivic Adams spectral sequence
\[
E_2 = \Ext_k \Rightarrow \pi^k_{*,*}.
\]
\index{motivic Adams spectral sequence}
\index{Adams spectral sequence, motivic}
\index{spectral sequence, motivic Adams}
We will not give a detailed construction of the spectral sequence.
Just as in the classical case, the idea is to construct an Adams resolution
of the motivic sphere spectrum, using copies of the 
motivic Eilenberg-Mac Lane $\MM \F_2$ and a few of its basic properties
\cite{Morel99} \cite{DI10} \cite{HKO11}.

A certain amount of care must be taken with 
convergence of the motivic Adams spectral sequence.
In particular, one must complete not only with respect to a prime,
but also with respect to the first Hopf map $\eta$.
Convergence results
are worked out in \cite{HKO11} \cite{HKO11b} \cite{Mantovani18}
\cite{KW18},
and we will not discuss them further here.

There are three major phases to carrying out an Adams spectral
sequence computation:
\begin{enumerate}
\item
Compute the $E_2$-page.  This is an algebraic problem.  
See Section \ref{subsctn:A-cohlgy}.
\item
Compute Adams differentials.  This is the hardest of the three steps.
\item
The $E_\infty$-page gives the filtration quotients in a filtration
on $\pi_{*,*}$, but these quotients can hide some of the structure
of $\pi_{*,*}$. 
We may have two elements
$\alpha$ and $\beta$ of $\pi_{*,*}$ detected by elements
$a$ and $b$ of the $E_\infty$-page with filtrations $i$ and $j$
respectively.
Sometimes, the product $\alpha \beta$ is non-zero in
$\pi_{*,*}$, but is detected in filtration greater than $i+j$.
In this case, the product $a b$ equals zero in the
$E_\infty$-page, even though $\alpha \beta$ is non-zero.
Some work is required to reconstruct
$\pi_{*,*}$ itself from its filtration quotients.
See \cite{Isaksen14c}*{Section 4.1.1} for a careful discussion
of these issues.
\end{enumerate}

\section{The motivic effective slice spectral sequence}
\label{sctn:tmsss}

In this section we review the effective slice filtration of the
stable motivic homotopy category $\SH(k)$.  
This filtration was originally
proposed by Voevodsky \cite{voevodsky.open} under the name
``slice filtration".
The goal is to understand the associated spectral sequence for the sphere. 
The latter approximates the stable motivic homotopy groups in a way we will make precise. 

The effective slice spectral sequence \eqref{equation:splicespectralsequenceforsphere} and its 
very effective 
version introduced in \cite{SO} are algebro-geometric analogues of the topological Atiyah-Hirzebruch 
spectral sequence \cite{zbMATH03176572}.
Among the applications of the effectice filtration we note Voevodsky's construction of the effective slice spectral sequence relating motivic cohomology to algebraic $K$-theory \cite{Voevodsky:motivicss}, 
Levine's proof of full faithfulness of the constant functor from the stable topological homotopy category to the stable motivic homotopy category over algebraically closed fields \cite{levine.comparison},
a new proof of Milnor's conjecture on quadratic forms in \cite{roendigs-oestvaer.hermitian}, 
\index{Milnor conjecture}
and the identification of Morel's plus part of the rational sphere spectrum with rational motivic cohomology along with finiteness for the motivic stable homotopy groups over finite fields in 
\cite{roendigs-spitzweck-oestvaer.1line}.

The equivariant slice filtration has become an important tool
classical homotopy theory.  For example, it is central to
major progress on elements of Kervaire invariant one \cite{HHR16}.
The reader should beware that 
the equivariant slice filtration does not directly correspond to
the motivic effective slice filtration.
However, the filtrations do happen to align in many cases 
related to $K$-theory and cobordism.

\subsection{The effective slice filtration}
\label{subsctn:effective-filtration}

Let $\SH^{\eff}(k)$ be the localizing subcategory of 
the stable motivic homotopy category 
$\SH(k)$ generated by the motivic suspension spectra 
of all smooth schemes over $k$,
i.e., 
the smallest full triangulated subcategory of $\SH(k)$ that contains suspension spectra of smooth schemes and is closed under coproducts \cite{voevodsky.open}.
These subcategories
form the effective slice filtration
\begin{equation}
\label{equation:slicefiltration}
\cdots
\subset
\Sigma^{2q+2,q+1}\SH^{\eff}(k)
{\subset}
\Sigma^{2q,q}\SH^{\eff}(k)
{\subset}
\Sigma^{2q-2,q-1}\SH^{\eff}(k)
\subset
\cdots.
\end{equation}
\index{effective slice filtration}
\index{slice filtration, effective}
Since $\SH^{\eff}(k)$ is closed under
simplicial suspsension and desuspension,
the subcategory 
$\Sigma^{2q,q}\SH^{\eff}(k)$ is equal to the subcategory
$\Sigma^{0,q}\SH^{\eff}(k)$.
The effective slice filtration is exhaustive in the sense that
$\SH(k)$ is the smallest triangulated subcategory that is
closed under coproducts and contains
$\Sigma^{2q,q}\SH^{\eff}(k)$ for all $q$.

The inclusion $\Sigma^{2q,q}\SH^{\eff}(k)\subset\SH(k)$ 
is left adjoint to
a functor that takes a motivic spectrum to its
``$q$th effective slice cover".  (This is analogous to the classical inclusion of
$q$-connected spectra into all spectra, which is left adjoint
to the $q$-connected cover.)
We write $\f_q$ for the functor that takes a motivic spectrum
to its $q$th effective slice cover.

Any motivic spectrum $\EE$ has an effective slice filtration
\begin{equation}
\label{eq:slice-filt}
\cdots \to \f_{q+1} \EE \to \f_q \EE \to \f_{q-1} \EE \to \cdots.
\end{equation}
The associated graded object of this filtration is
$\s_* \EE$, where the $q$th component $\s_q \EE$ is the cofiber
of
$\f_{q+1} \EE \to \f_q \EE$.

A related construction employs the very effective subcategory $\SH^{\veff}(k)$ introduced in \cite{SO}*{\S5}.
\index{very effective slice filtration}
\index{slice filtration, very effective}
This is the smallest full subcategory of $\SH(k)$ that contains all suspension spectra of smooth schemes of finite type, is closed under homotopy colimits, and is closed under extensions in the sense
that
if $X \to Y \to Z$ is a cofiber sequence such that 
$X$ and $Z$ belong to $\SH^{\veff}(k)$, then so does $Y$.

We note that $\SH^{\veff}(k)$ is contained in $\SH^{\eff}(k)$ but it is not a triangulated subcategory of $\SH(k)$ 
since it is not closed
under simplicial desuspension.
The very effective slice filtration takes the form
\begin{equation}
\label{equation:veffslicefiltration}
\cdots
\subset
\Sigma^{2q+2,q+1}\SH^{\eff}(k)
{\subset}
\Sigma^{2q,q}\SH^{\eff}(k)
{\subset}
\Sigma^{2q-2,q-1}\SH^{\eff}(k)
\subset
\cdots,
\end{equation}
and we obtain notions of very effective slice covers $\tilde{\f}_{q}$ and very effective
slices $\tilde{\s}_{q}$ for all $n\in\Z$.
To connect this to motivic homotopy groups we mention the fact that if $\EE\in\SH^{\veff}(k)$ then $\pi_{m,n}\EE=0$ for $m<n$ \cite{SO}*{Lemma 5.10}.
Some examples of very effective motivic spectra are algebraic cobordism $\MGL$, the effective slice cover $\kgl:=\f_{0}\KGL$ of algebraic $K$-theory \cite{SO},
and the very effective slice cover $\kq:=\tilde{\f}_{0}\KQ$ of hermitian $K$-theory \cite{ARO17}.
\index{algebraic cobordism spectrum}
\index{cobordism spectrum, algebraic}
\index{algebraic $K$-theory spectrum}
\index{$K$-theory spectrum, algebraic}
\index{hermitian $K$-theory spectrum}
\index{$K$-theory spectrum, hermitian}
One advantage of the very effective slice filtration is that it maps to the topological Postnikov tower under Betti realization \cite{GRSO}*{\S3.3}.

Higher structural properties of motivic spectra such as $A_{\infty}$- and $E_{\infty}$-structures interact well with these filtrations.
\index{$A_\infty$-structure}
\index{$E_\infty$-structure}
In \cite{GRSO} it is shown that if $\EE$ is an $A_{\infty}$- or $E_{\infty}$-algebra in motivic symmetric spectra \cite{Jardine} then $\f_{0}\EE$ and $\tilde{\f}_{0}\EE$ are naturally equipped 
with the structure of an $A_{\infty}$- resp. $E_{\infty}$-algebra.
Moreover, 
the canonical maps $\f_{0}\EE\to\EE$ and $\tilde{\f}_{0}\EE\to\EE$ can accordingly be modelled as maps of $A_{\infty}$- and $E_{\infty}$-algebras.
For every $q\in\Z$ the functors $\f_{q}$, $\s_{q}$, $\tilde{\f}_{q}$, and $\tilde{\s}_{q}$ respect module structures over $E_{\infty}$-algebras.

Applying motivic homotopy groups to the filtration
\eqref{eq:slice-filt}
yields by standard techniques an exact couple and the effective slice spectral sequence
\begin{equation}
\label{equation:splicespectralsequence}
E^{1}_{m,q,n}(\EE) 
=
\pi_{m,n}\s_{q}\EE
\Longrightarrow
\pi_{m,n}\EE.
\end{equation}
\index{effective slice spectral sequence}
\index{slice spectral sequence, effective}
\index{spectral sequence, effective slice}
For calculations it is important to note that the first
effective slice differential is induced by the composition $\s_{q}\EE\to\Sigma^{1,0}\f_{q+1}\EE\to\Sigma^{1,0}\s_{q+1}\EE$.

\subsection{The slices of $\unit$}
\label{subsctn:slice-1}

By \cite{roendigs-spitzweck-oestvaer.1line}, the effective slice spectral sequence for the sphere $\unit$ over any field has
$E_1$-page given by
\begin{equation}
\label{equation:splicespectralsequenceforsphere}
E^{1}_{m,q,n} = \pi_{m,n}\s_{q}(\unit).
\end{equation}
This spectral sequence converges conditionally, in the sense
of Boardman \cite{Boardman} to the $\eta$-completed stable motivic homotopy groups.
Technically, there is no need to complete at a prime, i.e.,
the effective slice spectral sequence is defined integrally.  However,
for practical purposes
it is best to further complete one prime at a time.
Here $\s_{q}(\unit)$ is the $q$th slice of $\unit$, i.e.,
the $q$th layer of the effective slice filtration on $\unit$.

The first step in understanding the effective slice spectral
sequence \eqref{equation:splicespectralsequenceforsphere} 
is to understand the slices $\s_q(\unit)$.
Voevodsky \cite{voevodsky.open}*{Conjectures 9 and 10}
predicted correctly that these slices are governed by
the structure of the Adams-Novikov $E_{2}$-page
(see \eqref{equation:sphereslices} below).
It turns out that $\s_{0}(\unit)$ identifies with the integral motivic Eilenberg-Mac Lane spectrum $\MM \Z$,
while $\s_{q}(\unit)$ for $q\geq 1$ are 
$\s_0(\unit)$-modules, i.e., $\MM \Z$-modules.
\index{motivic Eilenberg-Mac Lane spectrum}
\index{Eilenberg-Mac Lane spectrum, motivic}

We will now sketch how the slices $\s_q(\unit)$ can
be computed.
First, a geometric argument shows that 
the cone of the unit map $\unit\to\MGL$ lies in $\Sigma^{{2,1}}\SH^{\eff}(k)$.
It follows that $\s_{0}(\unit)\to \s_{0}(\MGL)$ is an isomorphism since $\s_{0}$ is trivial on $\Sigma^{{2,1}}\SH^{\eff}(k)$. 
\index{algebraic cobordism spectrum}
\index{cobordism spectrum, algebraic}

The slices of $\MGL$ are
described in Theorem \ref{thm:slice-MGL}.
It turns out to be more convenient to describe all the 
slices of $\MGL$ simultaneously, rather than one at a time.
So we will describe the entire graded motivic spectrum
$\s_* \MGL$.  We use the term ``effective slice degree" to refer to this external grading.

\begin{remark}
\label{remark:exponentialcharacteristic}
For the following discussion,  
if the base field $k$ is of positive characteristic,
one must invert its characteristic in the motivic spectra 
$\MGL$ and $\MM \Z$.
See \cite{roendigs-spitzweck-oestvaer.1line}*{\S2.1-2.2} for precise statements.
For legibility we will not make any notational changes.
\end{remark}

\begin{theorem}
\label{thm:slice-MGL}
\cite{spitzweck.relations}*{Corollary 4.7}
\cite{hoyois}*{\S8.3}
\cite{spitzweckalgebraiccobordism}*{Theorem 3.1}
Over a base field of characteristic zero the graded motivic spectrum $\s_* \MGL$ is equal to
$\MM \Z [x_1, x_2, \ldots]$,
where $x_j$ has motivic degree $(2j,j)$ and effective slice degree $j$.
\end{theorem}

The expression $\MM \Z [x_1, x_2, \ldots]$ in Theorem \ref{thm:slice-MGL}
bears some explanation.  This object is a direct sum of copies
of suspensions $\MM \Z$, indexed by monomials in the $x_j$.
Each monomial $x_1^{i_1} x_2^{i_2} \ldots x_r^{i_r}$ contributes
a copy of $\Sigma^{2q,q} \MM \Z$ to $\s_q \MGL$, where
$q = \sum_{j=1}^r j i_j$.
Note that $q$ is the degree of the monomial of the same name
in $\pi_* \MU = \Z[x_1, x_2, \ldots]$, 
where $\MU$ is the classical complex cobordism spectrum
and $x_j$ has degree $2j$ \cite{Milnor60}.
\index{complex cobordism spectrum}
\index{cobordism spectrum, complex}

The positive slices of $\unit$ are determined in several steps, starting with the standard cosimplicial $\MGL$-resolution 
\index{cosimplicial resolution}
\begin{displaymath}
\label{equation:hpoikn}
\xymatrix{
\unit
\ar[r] &
\MGL
\ar@<3pt>[r] \ar@<-3pt>[r] & 
\MGL^{\wedge 2} 
\ar@<6pt>[r] \ar[r] \ar@<-6pt>[r] & 
\MGL^{\wedge 3} \cdots, }
\end{displaymath}
which induces a natural isomorphism 
\begin{equation}
\label{equation:hdgwgv}
\s_{q}(\unit)
\overset{\simeq}{\to}
\underset{\Delta}{\holim}\,
\s_{q}(\MGL^{\wedge(\bullet+1)}).
\end{equation}

In order to describe the slices of
$\MGL^{\wedge n}$, first recall 
the smash product decompositions 
\begin{equation*}
\MGL \wedge \MGL = \MGL[b_1,b_2,b_3,\ldots]
\]
and
\[
\MU\wedge \MU = \MU[b_1,b_2,b_3,\ldots]
\end{equation*}
from \cite{NSO}*{Lemma 6.4(i)}
and 
\cite{Adams74}*{p.\ 87}
respectively. 
Here $b_j$ is the standard choice of generator in the dual Landweber-Novikov algebra \cite{Ravenel86}*{Theorem 4.1.11}, and 
the motivic degree of $b_j$ is $(2j,j)$, while the classical
degree of $b_j$ is $2j$.
These decompositions extend to higher smash powers in an 
obvious way.

From Theorem \ref{thm:slice-MGL} and these smash product decompositions,
it follows that the graded motivic spectrum
$\s_* ( \MGL^{\wedge n} )$ is
$\MM \Z \otimes \pi_* (\MU^{\wedge n})$ \cite{voevodsky.open}*{Conjecture 6}.
More precisely, $\pi_*(\MU^{\wedge n})$ is a free abelian group
on a set of generators, concentrated in even degrees.  
Each generator of $\pi_{2q}(\MU^{\wedge n})$ contributes a
summand $\Sigma^{2q,q} \MM \Z$ to
$\s_q ( \MGL^{\wedge n} )$.

With additional work detailed in \cite{roendigs-spitzweck-oestvaer.1line}*{\S2.2}, the above produces an isomorphism
\begin{equation}
\label{equation:sphereslices}
\s_* (\unit) \simeq
\MM \Z \otimes \Ext_{\MU_* \MU}^{*,*} (\MU_*, \MU_*).
\end{equation}
More precisely, each copy of $\Z/n$ in
$\Ext_{\MU_* \MU}^{p,2q} (\MU_*, \MU_*)$
contributes a summand $\Sigma^{2q-p,q} \MM \Z/n$ to
$\s_q (\unit)$.
Here $p$ is the Adams-Novikov filtration, i.e., the
homological degree of the $\Ext$ group, and
$2q$ is the total degree, i.e., the stem plus the 
Adams-Novikov filtration.

We note the slices are $\MM\Z$-modules and the effective slice $d_{1}$-differentials are maps between motivic Eilenberg-MacLane spectra.
These facts make the effective slice spectral sequence amenable to calculations over a base scheme affording an explicit description of its motivic cohomology ring along with the action of the 
motivic Steenrod algebra.

The $\Ext$ groups in \eqref{equation:sphereslices} 
form the $E_{2}$-page of the Adams-Novikov spectral sequence, 
which has been extensively studied by topologists \cite{Ravenel86}. 
\index{Adams-Novikov spectral sequence}
\index{spectral sequence, Adams-Novikov}
From an algebraic viewpoint, this is the cohomology of the ``Hopf algebroid" corepresenting the functor carrying a commutative ring to the groupoid of formal group laws over it 
and their strict isomorphisms \cite{Adams74}.
If $(s,t)\neq (0,0)$, then $\Ext_{\MU_\ast\MU}^{s,t}(\MU_\ast,\MU_\ast)$ is finite.

In practice, $\Ext_{\MU_* \MU}(\MU_*, \MU_*)$ is most
easily studied one prime at a time.
For the Brown-Peterson spectrum $\BP$ at a prime $p$, 
\index{Brown-Peterson spectrum}
we have 
$$
\Ext_{\MU_\ast \MU}^{s,t}(\MU_\ast,\MU_\ast)=\bigoplus_{p}\Ext_{\BP_\ast \BP}^{s,t}(\BP_\ast,\BP_\ast)
$$
when $(s,t) \neq (0,0)$,
so one should really study
$\Ext_{\BP_\ast \BP}^{s,t}(\BP_\ast,\BP_\ast)$ separately
for each prime.

This input provides a systematic way of keeping track of the direct summands of $\s_{q}(\unit)$.
For example, at $p = 2$,
the well-known elements $\alpha_1^q$ in
$\Ext^{q,2q}_{\BP_* \BP}(\BP_*, \BP_*)$
contribute $\Sigma^{q,q} \MM \Z/2$ summands to $\s_q(\unit)$,
and the element $\alpha_{2/2}$ in
\[
\Ext^{1,4}_{\BP_* \BP}(\BP_*, \BP_*)
\]
contributes
a $\Sigma^{3,2} \MM \Z/4$ summand to $\s_2(\unit)$.
At $p = 3$, the element
$\alpha_1$ in 
\[
\Ext^{1,4}_{\BP_* \BP}(\BP_*, \BP_*)
\]
contributes $\Sigma^{3,2} \MM \Z/3$ to $\s_2(\unit)$.
These copies of $\Sigma^{3,2} \MM \Z/4$ and
$\Sigma^{3,2} \MM \Z/3$ assemble into a copy of
$\Sigma^{3,2} \MM \Z/12$ in $\s_2(\unit)$.

\begin{remark}
\label{rmk:slice-mult}
Our notation for $\s_*(\unit)$ suggests that we are
describing $\s_*(\unit)$ as a ring object, 
but the notation is somewhat deceptive.
For example, consider the map 
$\s_1(\unit) \wedge \s_1(\unit) \to \s_2(\unit)$ at $p = 2$,
where 
$\s_1(\unit)$ is equivalent to $\Sigma^{1,1} \MM\Z/2$ corresponding
to $\alpha_1$, and $\s_2(\unit)$ is equivalent to
\[
\Sigma^{2,2} \MM\Z/2 \vee \Sigma^{3,2} \MM\Z/4
\]
corresponding to $\alpha_1^2$ and $\alpha_{2/2}$
respectively.

The source $\Sigma^{1,1} \MM\Z/2 \wedge \Sigma^{1,1} \MM\Z/2$
splits as a wedge of suspensions
of copies of $\MM\Z/2$, where the wedge is indexed by a basis
for the motivic Steenrod algebra.
Thus we want to calculate the map
\begin{equation}
\label{eq:slice-product-ex}
\Sigma^{2,2} \bigvee_{b \in B} \Sigma^{d_b} \MM\Z/2 \to 
\Sigma^{2,2} \MM\Z/2 \vee \Sigma^{3,2} \MM\Z/4,
\end{equation}
where $B$ is a basis for the motivic Steenrod algebra,
and $d_b$ is the bidegree of $b$.

If we restrict to the summand of the source of
\eqref{eq:slice-product-ex} corresponding to
the element $1$ of the motivic Steenrod algebra,
then we have a map
\[
\Sigma^{2,2} \MM\Z/2 \to 
\Sigma^{2,2} \MM\Z/2 \vee \Sigma^{3,2} \MM\Z/4.
\]
This map is the identity on the first summand of the target
and zero on the second summand of the target,
corresponding to the multiplicative relation
$\alpha_1 \cdot \alpha_1 = \alpha_1^2$.

However, if we restrict to the summand of the source 
of \eqref{eq:slice-product-ex} corresponding
to the element $\Sq^1$ of the motivic Steenrod algebra,
then we have a map
\[
\Sigma^{3,2} \MM\Z/2 \to 
\Sigma^{2,2} \MM\Z/2 \vee \Sigma^{3,2} \MM\Z/4.
\]
This map turns out to be non-trivial when restricted to the
second summand of the target, and this has explicit computational
consequences for the effective slice spectral sequence.

This warning about multiplicative structures applies
to other slice calculations in this manuscript.
\end{remark}

It is natural to ask how $\eta$ acts on the slices of $\unit$.
A complete description can be extracted from the work of Andrews-Miller \cite{AM14} that calculates the $\alpha_1$-periodic
$\Ext$ groups over $\BP_* \BP$ at the prime $2$.
This shows that $\s_* (\unit) [\eta^{-1}]$ is 
\begin{equation}
\label{equation:etainvertedsphereslices}
\MM \Z/2[\alpha_1^{\pm 1}, \alpha_3, \alpha_4]/\alpha_3^2,
\end{equation}
where $\alpha_1$ has motivic degree $(1,1)$ and effective slice degree $1$;
$\alpha_3$ has motivic degree $(5,3)$ and effective slice degree $3$;
and $\alpha_4$ has motivic degree $(7,4)$ and effective slice degree $4$.

\subsection{The slices of other motivic spectra}
\label{subsctn:slice-other}

Related to \eqref{equation:sphereslices}, 
the slices of algebraic $K$-theory $\KGL$, 
hermitian $K$-theory $\KQ$, 
and higher Witt-theory $\KW$ can be identified.

\begin{theorem}
\label{thm:slice-K}
\cite{Levine08}
\cite{roendigs-oestvaer.hermitian}
\begin{enumerate}
\item
\label{part:algebraicktheoryslices}
The graded motivic spectrum $\s_{*}\KGL$ is equivalent to
$\MM\Z[\beta^{\pm 1}]$, where $\beta$ has motivic degree $(2,1)$ 
and effective slice degree $1$.
\item
\label{part:hermitianktheoryslices}
When $\Char(k)\neq 2$, the
graded motivic spectrum $\s_* \KQ$ is equivalent to
\[
\MM\Z[\alpha_1, v_1^{\pm 2}] / 2 \alpha_1,
\]
where $\alpha_1$ has motivic degree $(1,1)$ and effective slice degree $1$,
and $v_1^2$ has motivic degree $(4,2)$ and effective slice degree $2$.
\item
\label{part:wittheoryslices}
When $\Char(k)\neq 2$, the
graded motivic spectrum $\s_* \KW$ is equivalent to
\[
\MM \Z/2 [\alpha_1^{\pm 1}, v_1^{\pm 2}],
\]
where $\alpha_1$ has motivic degree $(1,1)$ and effective slice degree $1$,
and $v_1^2$ has motivic degree $(4,2)$ and effective slice degree $2$.
\end{enumerate}
\end{theorem}
\index{algebraic $K$-theory spectrum}
\index{$K$-theory spectrum, algebraic}
\index{hermitian $K$-theory spectrum}
\index{$K$-theory spectrum, hermitian}
\index{higher Witt-theory spectrum}
\index{Witt-theory spectrum, higher}

Part (\ref{part:algebraicktheoryslices}) of Theorem
\ref{thm:slice-K} says that 
$\s_q \KGL$ is $\Sigma^{2q,q} \MM \Z$.

In Part (\ref{part:hermitianktheoryslices}) of Theorem
\ref{thm:slice-K},
the monomials $v_1^{2q}$ contribute copies of
$\Sigma^{4q, 2q} \MM \Z$ to $\s_{2q} \KQ$.
Also, for $p \geq 1$, 
monomials $\alpha_1^p v_1^{2q}$ contribute copies
of $\Sigma^{4q+p, 2q+p} \MM \Z/2$ to $\s_{2q+p} \KQ$.
Beware that the notation for $\s_* \KQ$ only partly
describes the multiplicative structure, as in 
Remark \ref{rmk:slice-mult}.

We obtain Part (\ref{part:wittheoryslices}) 
of Theorem \ref{thm:slice-K} from 
Part (\ref{part:hermitianktheoryslices})
since $\KW$ is obtained from $\KQ$ by inverting $\eta$
in the same way that
$\unit[\eta^{-1}]$ is obtained from $\unit$.

\begin{remark}
Spitzweck's work on motivic cohomology in \cite{Spitzweck} shows the isomorphisms in \eqref{equation:sphereslices}, 
\eqref{equation:etainvertedsphereslices}, 
and all three parts of Theorem \ref{thm:slice-K}
hold not only over fields, but also
over Dedekind domains of mixed characteristic with no residue fields of characteristic $2$
(see \cite{roendigs-spitzweck-oestvaer.1line}*{\S2.3}).
\end{remark}

Thanks to the slice calculations reviewed in this section we know precisely how the $E^{1}$-pages of the effective slice spectral sequences for $\unit$, $\unit[\eta^{-1}]$, $\KGL$, $\KQ$, 
and $\KW$ are given in terms of motivic cohomology groups of the base scheme.
Now all the fun can begin with determining the corresponding differentials and 
resolving the information hidden by
the associated graded structure of the
effective slice $E_\infty$-pages.

\section{$\C$-motivic stable homotopy groups}
\label{sctn:C}

In this section, we fix the base field $k = \C$.  
This special case is made easier by the fact 
that $\C$ has trivial arithmetic
properties, i.e., $K^M_*(\C)/p$ is trivial.  
In fact,
the calculations for 
any algebraically closed field of characteristic zero work out identically.

To begin, the stable homotopy groups $\pi_{*,0}^\C$ are isomorphic
to the classical stable homotopy groups $\pi_*$ \cite{levine.comparison}.
This isomorphism occurs even without completion.
However, complicated exotica phenomena occur in other weights.

In this section, we also fix the prime $p = 2$.  
There are interesting phenomena to study at odd primes, but
those cases have not yet been studied as extensively.
See \cite{Stahn16} for some $\C$-motivic results at odd primes.

As mentioned in Remark \ref{rem:A-cohlgy-C},
the cohomology of the $\C$-motivic Steenrod algebra is completely known
in a large range. 
Current computer calculations extend beyond the 100-stem, with
partial information out to the 140-stem.
We recommend that the reader view $\C$-motivic Adams
charts throughout this discussion \cite{Isaksen14a}.

\subsection{$\C$-motivic and classical stable
homotopy groups}
\label{subsctn:C-cl-compare}

As observed in Remark \ref{rem:A-invert-tau},
inverting $\tau$ takes the
$\C$-motivic Steenrod algebra to the classical Steenrod algebra
with $\tau^{\pm 1}$ adjoined.
In fact, this principle extends to $\Ext$ groups, so that there 
is an isomorphism
\[
\Ext_\C (\M_2, \M_2) [\tau^{-1}] \cong
\Ext_{\cl}(\F_2, \F_2) [\tau^{\pm 1}],
\]
where $\Ext_{\cl}(\F_2,\F_2)$ represents $\Ext$ groups
over the classical Steenrod algebra.
Even further, there is an isomorphism between the
$\tau$-periodic $\C$-motivic Adams spectral sequence and the
classical Adams spectral sequence.
Consequently, the $\tau$-periodic $\C$-motivic stable homotopy groups
are isomorphic to the classical stable homotopy groups
with $\tau^{\pm 1}$ adjoined.

This comparison between the $\C$-motivic and classical situations
is induced by the Betti realization functor that takes a 
complex variety to its underlying topological space of $\C$-valued points.

\subsection{$\eta$-periodic-phenomena}
\label{subsctn:C-eta-periodic}

Classically, the element $h_1^4$ of $\Ext_{\cl}$ is zero.
From Section \ref{subsctn:C-cl-compare}, it follows that
$\C$-motivically, $\tau^k h_1^4$ equals zero for some $k$.
Classically, there is a relation $h_0^2 h_2 = h_1^3$.
In the $\C$-motivic setting, 
the weight of $h_0^2 h_2$ is $2$, 
while the weight of $h_1^3$ is $3$.
Consequently, the correct motivic relation is
$h_0^2 h_2 = \tau h_1^3$.
Therefore, 
\[
\tau h_1^4 = h_0^2 h_1 h_2 = 0,
\]
but it turns out that $h_1^4$ is not zero.
Moreover, every element $h_1^k$ is non-zero for all $k \geq 0$.

There is a naive explanation for this phenomenon.
Recall that the element $h_1$ is detected by $[ \xi_1 ]$ in the
$\C$-motivic cobar resolution.  Unlike in the classical case,
the element $\xi_1$ is indecomposable in the $\C$-motivic 
dual Steenrod
algebra.
It is this property of the $\C$-motivic dual 
Steenrod algebra that ultimately
allows $h_1$ to be not nilpotent.

This exotic behavior of $h_1$ leads us inevitably to ask about the effect
of inverting $h_1$ on the cohomology of the $\C$-motivic Steenrod
algebra, and about the effect of inverting the first Hopf map $\eta$
in the $\C$-motivic stable homotopy groups.
Inspection of a motivic Adams chart \cite{Isaksen14a} reveals many
classes, such as $c_0$, $d_0$, and $e_0$ that are
non-zero after inverting $h_1$.

We draw particular attention to the element $B_1$ 
(also known as $M h_1$) in the $46$-stem.  
From the perspective of the May spectral sequence, it is a surprise that
this element survives $h_1$-inversion because it is unrelated to any of
the many elements in lower stems that survive $h_1$-inversion.  
This observation
eventually led to a clean calculation.

\begin{theorem}
\cite{GI15}
The $h_1$-periodic cohomology of the $\C$-motivic Steenrod algebra
is
\[
\Ext_\C [h_1^{-1}] \cong \F_2 [h_1^{\pm 1}] [ v_1^4, v_2, v_3, \cdots ],
\]
where $v_1^4$ has degree $(8,4,4)$ and $v_n$ has degree 
$(2^n-2, 1, 2^{n-1} - 1)$.
\end{theorem}

The article \cite{GI15} conjectured 
that there are Adams differentials $d_2(v_n) = h_1 v_{n-1}^2$ for all
$n \geq 3$.
This pattern of Adams differentials
would determine the $\eta$-periodic $\C$-motivic stable homotopy groups.
The conjecture was proved by Andrews and Miller \cite{AM14},
who computed the $\alpha_1$-periodic $E_2$-page of the
Adams-Novikov spectral sequence.
\index{Adams-Novikov spectral sequence}
\index{spectral sequence, Adams-Novikov}
This Adams-Novikov computation immediately determines
the $\eta$-periodic $\C$-motivic stable homotopy groups,
using the circle of ideas discussed below in 
Section \ref{subsctn:Adams-diff}.

\begin{theorem}
\label{thm:eta-C}
\cite{AM14}
The $\eta$-periodic $\C$-motivic stable homotopy groups are
\[
\F_2 [ \eta^{\pm 1} ] [ \mu, \sigma ] / \sigma^2,
\]
where $\mu$ has degree $(9,5)$ and $\sigma$ has degree $(7,4)$.
\end{theorem}

Theorem \ref{thm:eta-C} holds even without completions.

\subsection{Adams differentials}
\label{subsctn:Adams-diff}

The hardest part of a classical or motivic Adams spectral sequence computation is the
determination of Adams differentials. 
\index{motivic Adams spectral sequence}
\index{Adams spectral sequence, motivic}
\index{spectral sequence, motivic Adams}
Some techniques include:
\begin{enumerate}
\item
Relations obtained by shuffling Toda brackets,
\index{Toda bracket}
which then imply that differentials must occur.
\item
Moss's theorem \cite{Moss70}, which can show that certain elements
in the Adams spectral sequence 
must survive to detect certain Toda brackets.
\item
Bruner's theorem on the interaction between algebraic squaring operations
and Adams differentials \cite{Bruner84}.
\item
Comparison to other information about stable homotopy groups, such
as the image of $J$ \cite{Adams66}, the homotopy groups of
$\tmf$ \cite{tmf14}, and the Adams-Novikov spectral sequence.
\index{image of $J$}
\index{$J$, image of}
\index{$\tmf$}
\index{topological modular forms}
\index{Adams-Novikov spectral sequence}
\index{spectral sequence, Adams-Novikov}
\end{enumerate} 

Many examples of these techniques appear in \cite{BMT70} \cite{Bruner84}
\cite{Isaksen14c} \cite{IX15} \cite{MT67} \cite{WangXu16}.  
The manuscript \cite{Isaksen14c}
contains an exhaustive discussion of Adams differentials up to the
59-stem.  With much effort, 
this bound was pushed slightly further to the 61-stem \cite{WangXu16}.
However, all of these techniques become impractical
in higher stems.

Recent work on the comparison between the $\C$-motivic and classical
stable homotopy categories has provided a new tool for computing
Adams differentials that has allowed us to extend
computations into a much larger range \cite{IWX18}.
Current calculations extend beyond the 90-stem,
with ongoing progress into even higher stems.

Two new ingredients allow for this program to succeed.
The first is computer data for the Adams-Novikov $E_2$-page.
The second ingredient is the following theorem.

\begin{theorem}
\label{thm:Ctau}
\cite{Gheorghe17b}
\cite{GWX18}
\cite{Pstragowski18}
In $\C$-motivic stable homotopy theory,
\begin{enumerate}
\item
The cofiber $\unit/\tau$ of $\tau$ is an $E_\infty$-ring spectrum,
in an essentially unique way.
\index{cofiber of $\tau$}
\index{$\tau$, cofiber of}
\item
With suitable finiteness conditions, the homotopy category
of $\unit/\tau$-modules is equivalent to the derived category
of $\BP_* \BP$-comodules.
\item
The Adams spectral sequence for $\unit/\tau$ is isomorphic to
the algebraic Novikov spectral sequence 
\cite{Novikov67} \cite{Miller75}
that converges to the Adams-Novikov $E_2$-page.
\index{algebraic Novikov spectral sequence}
\index{Novikov spectral sequence, algebraic}
\index{spectral sequence, algebraic Novikov}
\end{enumerate}
\end{theorem}

It is remarkable and unexpected that a homotopical category 
such as $\unit/\tau$-modules would have a purely algebraic 
description.
The proof of part (3) of Theorem \ref{thm:Ctau} is that the
filtrations associated to the two spectral sequences correspond
under the equivalence of part (2).

A new, more powerful approach to computing Adams differentials
is summarized in the following steps:
\begin{enumerate}
\item
Compute
the cohomology of the $\C$-motivic Steenrod algebra by machine.  
\item
Compute by machine
the algebraic Novikov spectral sequence, including
all differentials and multiplicative structure.
\item
Use Theorem \ref{thm:Ctau}
to completely describe
the motivic Adams spectral sequence for $\unit/\tau$.
\item
Use the cofiber sequence
\[
\Sigma^{0,-1} \unit \stackrel{\tau}{\to} \unit \to \unit/\tau \to 
\Sigma^{1,-1} \unit
\]
and naturality of Adams spectral sequences 
to pull back and push forward Adams differentials for $\unit/\tau$
to Adams differentials for the motivic sphere.
\item
\label{item:ad-hoc}
Apply a variety of ad hoc arguments to deduce additional Adams differentials
for the motivic sphere, as described at the beginning of
this section.
\item
Use a long exact sequence in homotopy groups to deduce hidden $\tau$
extensions in the motivic Adams spectral sequence for the sphere.
\item
Invert $\tau$ to obtain the classical Adams spectral sequence and the 
classical stable homotopy groups.
\end{enumerate}

The weak link in this algorithm is
step (\ref{item:ad-hoc}), where ad hoc arguments come into play.
The method will continue to calculate new stable
homotopy groups until the ad hoc arguments become too 
complicated to resolve.
It is not yet clear when this will occur.

\subsection{Motivic nilpotence and periodicity}

Section \ref{subsctn:C-eta-periodic} 
discussed the non-nil\-pot\-ence of the 
element $\eta$ of $\pi_{1,1}^\C$.
This non-nilpotent behavior led Haynes Miller to propose that there
might be an infinite family of periodicity operators $w_n$ such that
$w_0$ is $\eta$, in analogy to the $v_n$-periodicity operators that
begin with the non-nilpotent element $v_0 = 2$.
\index{motivic periodicity}
\index{periodicity, motivic}
\index{motivic nilpotence}
\index{nilpotence, motivic}

In fact, Miller's guess turned out to be correct \cite{Andrews14}
\cite{Gheorghe17}.  In more detail, there exist $\C$-motivic ring
spectra $\K(w_n)$ whose motivic stable homotopy groups are of the form
$\F_2[w_n^{\pm 1}]$.  Also, for each $n$, 
there exist finite $\C$-motivic complexes
that possess $w_n$-self maps.  More precisely, there exist
complexes $\mathbf{X}_n$ equipped with maps of the form
\[
\xymatrix@1{
\Sigma^{d |w_n|} \mathbf{X}_n \ar[r]^-{f_n} & \mathbf{X}_n 
}
\]
that induce isomorphisms in $\K(w_n)$-homology.
One can then study $w_n$-periodic families of elements in
$\pi_{*,*}^\C$ by considering compositions of the form
\[
\xymatrix@1{
\Sigma^{p,q} \unit \ar[r] & 
\Sigma^{d m |w_n|} \mathbf{X}_n \ar[r]^{f_n} & \cdots \ar[r]^{f_n} &
 \mathbf{X}_n \ar[r] &
\unit.}
\]
Andrews \cite{Andrews14} found the first explicit examples
of $w_1$-periodic families.

It turns out that 
there are even more exotic 
$\C$-motivic periodicities \cite{Krause18}.
Roughly speaking, these periodicities correspond to the elements
$h_{ij}$ of the May spectral sequence for $i \geq j$.
The $v_n$-periodicities correspond to the elements $h_{n+1,0}$,
while the $w_n$-periodicities correspond to
$h_{n+1,1}$.
While there do exist $\C$-motivic spectra of the form $\K(h_{ij})$,
these objects do not all possess ring structures.

\section{$\R$-motivic stable homotopy groups}
\label{sctn:R}

In this section, we fix the base field $k = \R$.
This special case introduces new phenomena not seen in the
$\C$-motivic situation.  The first complication is
that
$K^M_*(\R)/2$ is now a polynomial algebra on one class $\rho$.
Studying the $\R$-motivic case allows us to grapple with
the difficulties presented by the non-zero element $\rho$.
However, it avoids even further complications 
that would be created by relations involving $\rho$.

In this section, we also fix the prime $p = 2$.
The $\R$-motivic case at odd primes is not so interesting, since
$K^M_*(\R)/p$ is zero.

As mentioned in Section \ref{subsctn:rho-Bockstein}, the 
$\rho$-Bockstein spectral sequence is an effective tool
for computing the cohomology of the 
$\R$-motivic Steenrod algebra.
Preliminary work shows that the spectral sequence can be
completely analyzed in a large range.
See \cite{DI17} for a sense of the structure
of the calculation, although the range studied there is just
a tiny portion of what is possible.
\index{$\rho$-Bockstein spectral sequence}
\index{spectral sequence, $\rho$-Bockstein}
\index{Bockstein spectral sequence, $\rho$}

\subsection{$\eta$-periodic phenomena}
\label{subsctn:R-eta-periodic}

Just as in the $\C$-motivic case,
one can study the effect of inverting $h_1$ 
on the cohomology of the $\R$-motivic Steenrod algebra,
and about the effect of inverting the first Hopf map $\eta$
in the $\R$-motivic stable homotopy groups.
This study was carried out to completion in \cite{GI16}.

The first step is to consider the $h_1$-periodic version
of the $\rho$-Bockstein spectral sequence discussed in
Section \ref{subsctn:rho-Bockstein}, which takes the form
\[
\Ext_\C [ \rho ][h_1^{-1}] \Rightarrow \Ext_\R [h_1^{-1}].
\]
This spectral sequence can be completely analyzed.

The next step is to consider the $h_1$-periodic version
of the $\R$-motivic Adams spectral sequence, which takes
the form
\[
\Ext_\R [h_1^{-1}] \Rightarrow \pi_{*,*}^\R[\eta^{-1}].
\]
Once again, this spectral sequence can be completely analyzed,
but it is much more interesting than
the $\C$-motivic case discussed in 
Section \ref{subsctn:C-eta-periodic},
which collapses at the $E_3$-page.
In fact, in this case there are non-trivial differentials
on every page of the spectral sequence.  
These differentials are deduced by analyzing higher
homotopical structure, i.e., Toda brackets.

Before we can state the conclusion of the calculation,
we need some additional notation.  This complication arises
because we want to invert an element $\eta$ that has
non-zero degree $(1,1)$.
We write
$\Pi_m^\R$ for the direct sum 
\[
\bigoplus_{n \in \Z} \pi_{m+n, n}^\R,
\]
which we call the $m$th Milnor-Witt stem.
Then $\eta$ acts on each $\Pi_m^\R$, and we can consider
$\Pi_m^\R [\eta^{-1}]$.

\begin{theorem}
\label{thm:R-eta-periodic}
\cite{GI16}
The $\eta$-periodic $\R$-motivic stable homotopy groups
$\pi_{*,*}[\eta^{-1}]$ are given by:
\begin{enumerate}
\item
$\Pi_0^\R[\eta^{-1}] = \Z_2 [\eta^{\pm 1}]$.
\item
$\Pi_{4m-1}^\R[\eta^{-1}] = \Z/2^{u+1}[\eta^{\pm 1}]$
for $m > 1$, where $u$ is the $2$-adic valuation of $4m$.
\item
$\Pi_m^\R[\eta^{-1}] = 0$ otherwise.
\end{enumerate}
\end{theorem}

The groups that appear in Theorem \ref{thm:R-eta-periodic} 
are reminiscent of the groups that appear in the classical
image of $J$.  One might speculate on a more direct connection
between these $\eta$-periodic $\R$-motivic stable homotopy groups
and the image of $J$.  However, the Toda bracket structures
are different.
\index{image of $J$}
\index{$J$, image of}

\section{Motivic stable homotopy groups over general fields}
\label{sctn:general-field}

Sections \ref{sctn:C} and \ref{sctn:R} discussed what is 
known about $\C$-motivic and $\R$-motivic stable homotopy groups.
We now consider motivic stable homotopy groups over larger
classes of fields.  Naturally, specific information is harder
to obtain when the base field is allowed to vary widely.

In addition to the phenomena described in the following sections,
we also mention the article \cite{DI13}, which discusses
the existence of motivic Hopf maps and some relations satisfied
by these maps.

\subsection{Motivic graded commutativity}

In the classical case, the stable homotopy groups $\pi_*$ are 
graded commutative, in the sense that 
\[
\alpha \beta = (-1)^{|\alpha||\beta|} \beta \alpha
\]
for all $\alpha$ and $\beta$.  In other words,
$\alpha$ and $\beta$ anti-commute if they are both odd-dimensional
classes, and they strictly commute otherwise.
Ultimately, the graded commutativity arises from the fact that
the twist map
\[
S^1 \smash S^1 \to S^1 \smash S^1
\]
has degree $-1$.

\begin{remark}
\label{rem:2sigma^2}
The graded commutativity of the classical stable homotopy groups
has non-trivial consequences for the structure of the
Adams spectral sequence.  For example, consider the
third Hopf map $\sigma$ in $\pi_7$, which is detected by $h_3$
in the Adams spectral sequence.
Graded commutativity implies that $2 \sigma^2$ must equal
zero in $\pi_{14}$.  Now $h_0 h_3^2$ detects $2 \sigma^2$,
and it is non-zero in the Adams $E_2$-page.  Therefore,
there must be an Adams differential $d_2(h_4) = h_0 h_3^2$.
This is the first differential in the Adams spectral sequence.
\end{remark}

In the motivic situation, graded commutativity takes a slightly
more complicated form.
The twist map
\[
S^{1,0} \smash S^{1,0} \to S^{1,0} \smash S^{1,0}
\]
represents $-1$ in $\pi_{0,0}^k$, but the twist map
\[
S^{1,1} \smash S^{1,1} \to S^{1,1} \smash S^{1,1}
\]
represents a different element, usually called $\epsilon$,
in $\pi_{0,0}^k$.
This description of the twist maps leads to the following
theorem on motivic graded commutativity.
\index{motivic graded commutativity}
\index{graded commutativity, motivic}
\index{commutativity, motivic graded}

\begin{theorem}
\cite{Dugger14}*{Proposition 1.18}
\cite{DI13}*{Proposition 2.5}
\cite{Morel04}*{Corollary 6.1.2}
\label{thm:epsilon-commute}
Let $k$ be any base field.
For $\alpha$ in $\pi_{a,b}^k$ and $\beta$ in $\pi_{c,d}^k$,
\[
\alpha \beta = (-1)^{(a-b)(c-d)} \epsilon^{bd} \beta \alpha.
\]
\end{theorem}

\begin{remark}
Classically, the element $2$ of $\pi_0$ plays the role of the
``zeroth Hopf map".  Motivically, it is $1 - \epsilon$ that
plays this role.
The cofiber of $1-\epsilon$ is a 2-cell complex.  In the cohomology
of this 2-cell complex, there is a $\Sq^1$ operation that
connects the two cells.
In the cofiber of $2$, there is a $\Sq^1$ operation, but there
is also a non-trivial $\Sq^2$ operation when $\rho$ is non-zero.
\end{remark}

\subsection{Milnor-Witt $K$-theory}
\label{subsctn:Milnor-Witt}

In this section, we will recall the work of Morel on 
the motivic stable homotopy groups $\pi_{n,n}^k$
over arbitrary base fields $k$.
For more details,
see \cite{Morel04} \cite{Morel04b} \cite{Morel12}.
For related work, see also
\cite{DK18}
\cite{Hornbostel18}*{Appendix} \cite{Neshitov18} \cite{Thornton16}.

\begin{definition}
\label{defn:MWK}
\index{Milnor-Witt $K$-theory}
\index{$K$-theory, Milnor-Witt}
The Milnor-Witt $K$-theory 
$K^{MW}_*(k)$ of a field $k$ is the graded associative ring
generated by elements $[u]$ for all $u$ in $k^\times$, and
the element $\eta$, subject to the relations
\begin{enumerate}
\item
$[u] [1-u] = 0$ for all $u$ in $k$ except for $0$ and $1$.
\item
$[uv] = [u] + [v] + \eta [u] [v]$ for all $u$ and $v$ in $k^\times$.
\item
$[u] \eta = \eta [u]$ for all $u$ in $k^\times$.
\item
$\eta ( \eta[-1] + 2) = 0$.
\end{enumerate}
The degree of $[u]$ is $1$, and the degree of $\eta$ is $-1$.
\end{definition}

Setting $\eta = 0$ in Milnor-Witt $K$-theory recovers 
ordinary Milnor $K$-theory (see Section \ref{subsctn:Milnor-K}),
\index{Milnor $K$-theory}
\index{$K$-theory, Milnor}
so we can view
Milnor-Witt $K$-theory as a kind of deformation of Milnor
$K$-theory.
The first relation in Definition \ref{defn:MWK} is precisely
the Steinberg relation, 
\index{Steinberg relation}
while the second relation
is analogous to the additivity relation in Milnor $K$-theory,
with an error term involving $\eta$.

\begin{theorem}
\label{thm:Morel-MW0}
\cite{Morel04}*{Theorem 6.4.1}
\cite{Morel04b}*{Theorem 6.2.1}
\cite{Morel12}*{Corollary 1.25}
For any field $k$, the motivic stable homotopy group
$\pi_{n,n}^k$ is isomorphic to the $n$th Milnor-Witt
$K$-theory group $K^{MW}_n(k)$.
\end{theorem}

Note, in particular, that
$\pi_{0,0}^k$ is isomorphic to the Grothendieck-Witt ring
$GW(k)$
of quadratic forms over $k$.
Unlike most of the other results in this article,
there is no need for completions in the statement
of Theorem \ref{thm:Morel-MW0}.

The element $\epsilon$ that governs graded
commutativity corresponds here to $-1 - \eta [-1]$,
so the last relation of Definition \ref{defn:MWK}
is equivalent to the relation $\eta (1 - \epsilon) = 0$.
The element $\rho$ in $\pi_{-1,-1}$ corresponds to
$[-1]$, as in Section \ref{subsctn:Milnor-K}.

\subsection{The first Milnor-Witt stem}
\label{subsctn:1-line}

For simplicity let $k$ be a field of characteristic zero.
By the work reviewed in Section \ref{subsctn:Milnor-Witt} a next logical step is to compute the first Milnor-Witt stem
$\Pi^k_1 = \bigoplus_{n\in\Z}\pi_{n+1,n}^k$.
One of the major inspirations for this calculation is Morel's $\pi_{1}$-conjecture in weight zero, which states that
there is a short exact sequence
\begin{equation}
\label{eq:first-stable-stem2} 
0 
\to 
K^{M}_{2}(k)/24 
\to 
\pi_{1,0}^k
\to 
k^{\times}/2\oplus\Z/2
\to
0.
\end{equation}
This conjecture is proved in \cite{roendigs-spitzweck-oestvaer.1line}.

The kernel 
$K^{M}_{2}(k)/24$
in (\ref{eq:first-stable-stem2}) is
generated by the second motivic Hopf map $\nu$ in $\pi_{3,2}^k$, in the sense
that its elements are all of the form $\alpha \nu$, where
$\alpha$ is an element of $\pi_{-2,-2}^k$.
Such elements $\alpha$ correspond to elements of $K^M_2(k)$,
as in Theorem \ref{thm:Morel-MW0}.

On the other hand,
the image $k^{\times}/2\oplus\Z/2$ 
has two generators.
The second factor is generated by $\eta_{\Top}$, i.e.,
the image of the classical first Hopf map in $\pi^k_{1,0}$.
The first factor is generated by $\eta \eta_{\Top}$, in the sense
that its elements are all of the form $\alpha \eta \eta_{\Top}$,
where $\alpha$ is an element of $\pi_{-1,-1}^k$.
These generators are subject to the relations $24\nu=0$ and 
$12\nu=\eta^{2}\eta_{\Top}$, which are
related to the corresponding classical
relations $24 \nu = 0$ and $12 \nu = \eta^3$ in $\pi_3$.

It turns out the surjection in \eqref{eq:first-stable-stem2} arises from the unit map for the hermitian $K$-theory spectrum $\KQ$ of quadratic forms.
\index{hermitian $K$-theory spectrum}
One may speculate that its kernel is the image of a motivic $J$-homomorphism $K^{M}_{2}(k)=\pi_{1,0} GL\to\pi_{1,0}\unit$ for the general linear group $GL$. 
The Hopf construction should witness that $\nu$ is in the image of a motivic $J$-homomorphism, 
\index{image of $J$}
\index{$J$, image of}
so the relation $24\nu=0$ may be a shadow of some motivic version of the Adams conjecture. 
\index{Adams conjecture}

More generally, 
for every $n\in\Z$,
there is an exact sequence 
\cite{roendigs-spitzweck-oestvaer.1line}
\begin{equation}
\label{eq:first-stable-stem} 
0 
\to 
{K}^{M}_{2-n}/24(k)
\to 
\pi_{n+1,n}^k
\to 
\pi_{n+1,n}\f_{0}\KQ.
\end{equation}
Here $\f_{0}\KQ$ is the effective slice cover of hermitian $K$-theory 
(see Section \ref{sctn:tmsss}).
Note that the homotopy groups of $\KQ$ and $\f_{0}\KQ$ agree in nonnegative weight.
The rightmost map in \eqref{eq:first-stable-stem} is surjective for $n\geq -4$. 
In fact, the rightmost map is surjective for all $n$ if
$\f_0 \KQ$ is replaced by the very effective slice cover
$\kq = \tilde{\f}_0 \KQ$ \cite{RSO-2line}.
\index{very effective slice filtration}
\index{slice filtration, very effective}

The proof of \eqref{eq:first-stable-stem} is achieved by performing calculations with the effective slice spectral sequence for the sphere $\unit$, 
converging conditionally to the homotopy of the $\eta$-completion ${\unit}^\wedge_\eta$.
For all integers $n$,
there is a canonically induced isomorphism
\[ 
\pi_{n+1,n}\unit
\xrightarrow{\cong}
\pi_{n+1,n}{\unit}^\wedge_\eta
\]
noted in \cite{roendigs-spitzweck-oestvaer.1line}*{Corollary 5.2}.

The exact sequence \eqref{eq:first-stable-stem} vastly generalizes computations in \cite{OO14} for fields of cohomological dimension at most two,
i.e., 
if Milnor $K$-theory is concentrated in degrees $0$, $1$, and $2$.  
Examples include algebraically closed fields, finite fields, $p$-adic fields, and totally imaginary number fields (but not $\Q$ or $\R$).

It is interesting to compare the above with the computations of unstable motivic homotopy groups of punctured affine spaces in \cite{AF15} and \cite{asok-wickelgren-williams}.
If $d>3$, 
the extension for $\pi_d(\mathbf{A}^{d}\smallsetminus\{0\})$ conjectured by Asok-Fasel~\cite{mfo2013}*{Conjecture~7, p.~1894} coincides with \eqref{eq:first-stable-stem}. 
As noted in \cite{AF15}, 
the sheaf version of the exact sequence \eqref{eq:first-stable-stem} and a conjectural Freudenthal $\PP^1$-suspension theorem imply Murthy's conjecture on splittings of 
vector bundles \cite{AF15}*{Conjecture 1}.

\subsection{Finite fields}

This section describes work of Wilson and \O stv\ae r
\cite{WO17} on stable motivic homotopy groups over
finite fields of order $q = p^n$, 
where $p$ is an odd prime.

As always, the starting point is the motivic cohomology of a point.
Using Example \ref{ex:KM-Fq} and Theorem \ref{thm:cohlgy-point}, 
we obtain that
$H^{*,*}(\F_q; \F_2) = \F_2[u]/ u^2$.

At this point, the discussion splits into two cases.
If $p \equiv 1 \bmod 4$,
then $-1$ is a square in $\F_q$, and $\rho$ equals $0$.
Consequently, $\Ext$ calculations (i.e., motivic Adams
$E_2$-pages) are essentially the same as in the
$\C$-motivic case discussed in Section \ref{subsctn:A-cohlgy}.

On the other hand,
if $p \equiv 3 \bmod 4$,
then $u$ can be taken to be $\rho$.
The $\Ext$ groups in this case can be computed with
the $\rho$-Bockstein spectral sequence as in 
Section \ref{subsctn:rho-Bockstein}.
\index{$\rho$-Bockstein spectral sequence}
\index{spectral sequence, $\rho$-Bockstein}
\index{Bockstein spectral sequence, $\rho$}

However, there are significant differences to the
$\R$-motivic case that arise from the relation $\rho^2 = 0$.

In either case, there are interesting motivic Adams
differentials for finite fields that have no 
$\C$-motivic nor $\R$-motivic analogue
\cite{WO17}*{Corollary 7.12}.
More specifically, there are differentials of the form
\[
d_r (\tau^k) = u \tau^{k-1} h_0^r
\]
for some values of $r$ and $k$ that depend on the order
$q$ of the finite field $\F_q$.
These differentials are remarkable because they occur
already at the very beginning of the
spectral sequence in the $0$-stem!
The proofs of these differentials depend on a priori
knowledge of the motivic cohomology of $\F_q$ via its \'etale cohomology
\cite{Soule79}, as discussed in Section \ref{subsctn:MEMs}.
\index{etale cohomology}
\index{cohomology, etale}

An analysis of the 
Adams spectral sequence leads to an isomorphism
$$
\pi^{\F_{q}}_{n,0} 
\cong 
\pi_n^s \oplus\pi_{n+1}^s.
$$
for $0\leq n\leq 18$.  
In particular, the groups $\pi^{\F_{q}}_{4,0}$ and $\pi^{\F_{q}}_{12,0}$ are trivial.
It is interesting to note that the pattern $\pi^{\F_{q}}_{n,0} \cong \pi_n^s \oplus\pi_{n+1}^s$ does not hold in general.  
In fact when $q \equiv 5 \bmod 8$ one obtains
$$
\pi^{\F_{q}}_{19,0} \cong (\Z/8 \oplus \Z/2) \oplus \Z/4
\text{ and } \pi^{\F_{q}}_{20,0} \cong \Z/8 \oplus \Z/2.
$$

\subsection{$\eta$-periodic phenomena}

Relatively little is known about $\eta$-periodic
phenomena over fields other than $\C$ and $\R$.
Sections \ref{subsctn:C-eta-periodic}
and \ref{subsctn:R-eta-periodic} describe the
$\eta$-periodic groups
$\pi_{*,*}^{\C} [\eta^{-1}]$ and
$\pi_{*,*}^{\R} [\eta^{-1}]$.
Recent work of Wilson \cite{Wilson17}
calculates
$\pi_{*,*}^{k} [\eta^{-1}]$ 
for $\Q$ and for 
fields of cohomological dimension
at most 2, i.e., for fields whose Milnor $K$-theory vanishes
above degree 2.
All finite fields satisfy this condition.
The $\eta$-periodic groups are described in terms of the
Witt group $W(k)$ of quadratic forms over $k$.

For more general fields, 
the following theorem summarizes what we do know.
Recall the notation
\[
\Pi_m^k = \bigoplus_{n \in \Z} \pi_{m+n, n}^k
\]
from Section \ref{subsctn:R-eta-periodic}.

\begin{theorem}
\label{eta-periodic-rational}
\cite{ALP17}
Let $k$ be a field such that $\Char k \neq 2$.
Then 
\[
\Pi_m^k[\eta^{-1}] \otimes \Q = 0
\]
for all $m > 0$.
\end{theorem}

Theorem \ref{eta-periodic-rational}
leaves open the question of torsion in $\Pi_m^k[\eta^{-1}]$.
This torsion is likely to be quite interesting.

Additionally, we have some low-dimensional information.

\begin{theorem}
\cite{Roendigs-Tata}
Let $k$ be a field such that $\Char k \neq 2$.
Then $\Pi_1^k$ and $\Pi_2^k$ are both zero.
\end{theorem}

\section{Other motivic spectra}
\label{sctn:other}

One way to obtain information about stable motivic homotopy groups
is to consider other motivic spectra that are equipped with unit
maps from the motivic sphere spectrum.
The homotopy groups of these other
motivic spectra can give information about the motivic
stable homotopy groups by passing along the unit map.
We will discuss a few examples in this section.

Each motivic spectrum discussed below is of fundamental interest in
its own right.
Their associated cohomology theories detect interesting
phenomena in algebraic geometry, but that is not the focus of this
discussion.
We will not discuss their construction and geometric origins
because these topics go beyond
the scope of this article.

The unit map for algebraic cobordism factors through $\unit/\eta$, 
so every module over an oriented motivic ring spectrum is $\eta$-complete \cite{rso.hlp}*{Lemma 2.1}.
Hence the motivic Eilenberg-Mac Lane spectrum $\MM \Z$, 
the $K$-theory spectrum $\KGL$ and its covers, 
and the truncated Brown-Peterson spectra $\MBP\langle{n}\rangle$ are all $\eta$-complete.
On the other hand, 
hermitian $K$-theory $\KQ$ and higher Witt-theory $\KW$ do not coincide with their respective $\eta$-completions.

\subsection{Motivic Eilenberg-Mac Lane spectra}
\label{subsctn:MEMs}

Suppose $\ell$ is prime to the characteristic of the base field $k$.
Our first example is the motivic Eilenberg-Mac Lane spectrum $\MM \Z/{\ell^{\nu}}$ that represents motivic cohomology with $\Z/{\ell^{\nu}}$ coefficients \cite{Voevodsky98} \cite{RO08}.
\index{motivic Eilenberg-Mac Lane spectrum}
\index{Eilenberg-Mac Lane spectrum, motivic}
The stable homotopy groups of $\MM\Z/{\ell^{\nu}}$ are the same as the cohomology of a point with coefficients in $\Z/{\ell^\nu}$,
as described in Section \ref{subsctn:Milnor-K}.  
The mod-$\ell$ motivic Steenrod algebra is the ring of operations on $\MM\Z/{\ell}$, 
as described in Section \ref{subsctn:A-cohlgy}.  
This circle of ideas leads eventually to the motivic Adams spectral sequence of Section \ref{subsctn:mot-Adams}.

The Beilinson-Lichtenbaum conjecture, which is a consequence of the Milnor and Bloch-Kato conjectures, offers the following comparison isomorphism due to Voevodsky.
\index{Beilinson-Lichtenbaum conjecture}
\index{Bloch-Kato conjecture}
\index{Milnor conjecture}

\begin{theorem}
\label{thm:mot-et}
\cite{Voevodsky11}
For $p \leq q$ and $X$ a smooth $k$-scheme,
the \'{e}tale sheafification functor induces an isomorphism 
\begin{equation} 
H^{p,q}(X;\Z/\ell^{\nu}) 
\overset{\cong}{\rightarrow} 
H^p_{\text{\'et}}(X; \mu_{\ell^{\nu}}^{\otimes q})
\end{equation} 
between motivic and \'{e}tale cohomology groups,
where $\mu_{\ell^\nu}$ is the sheaf of $\ell^\nu$th roots of unity.
\end{theorem}
\index{etale cohomology}
\index{cohomology, etale}

This important isomorphism identifies the mod-$\ell^{\nu}$ motivic cohomology of $k$ with the classical cohomology groups of the absolute Galois group of $k$ \cite{MR1867431}.

Let 
$e$ be the exponent of
the multiplicative group $(\Z/\ell^{\nu})^{\times}$ of units.
The sheaf
$\mu_{\ell^\nu}^{\otimes e}$ is constant,
so
$H^{0,e}(k;\Z/\ell^\nu)$ is isomorphic to
$H^{0,0}(k; \Z/\ell^\nu) = \Z/\ell^\nu$.
Let us choose a generator $\tau_{\ell^{\nu}}$ of 
$H^{0,e}(k;\Z/\ell^{\nu})$.
Using Theorem \ref{thm:mot-et}, 
we deduce an isomorphism
\begin{equation} 
\label{equation:BLbottinverted}
H^{p,q}(X; \Z/\ell^{\nu})[\tau_{\ell^{\nu}}^{-1}] 
\cong 
H^p_{\text{\'et}}(X;\mu_{\ell^{\nu}}^{\otimes q})
\end{equation} 
for all integers $p,q\in\Z$ and all smooth $k$-schemes $X$.
The periodicity discussed here is related to
Thomason's seminar paper \cite{MR826102}, as well as
to a recent generalization to motivic spectra such as algebraic cobordism
\cite{elmanto-levine-spitzweck-oestvaer}.

In general, the motivic cohomology groups
$H^{p,q}(k; \Z/\ell^\nu)$ do not have a simple description
in terms of Milnor $K$-theory.
The case $\nu = 1$ is discussed in 
Section \ref{subsctn:Milnor-K}.
If $k$ contains all $\ell^\nu$th roots of unity,
then the sheaf 
$\mu_{\ell^{\nu}}$
is trivial, 
and Theorem \ref{thm:mot-et} gives a practial way of computing
$H^{p,q}(k; \Z/\ell^\nu)$.
One way to view these difficulties is that the
$\ell$-Bockstein spectral sequence that starts
with $H^{*,*}(k; \Z/\ell)$ and converges to
$H^{*,*}(k; \Z/\ell^\nu)$ is non-trivial.
\index{Bockstein spectral sequence}
\index{spectral sequence, Bockstein}

In addition to $\MM \Z/{\ell}$, we may also consider the motivic Eilenberg-Mac Lane spectrum $\MM \Z$ that represents integral motivic cohomology.
The cofiber sequence
\[
\xymatrix@1{
\MM \Z \ar[r]^{\ell^{\nu}} & \MM \Z \ar[r] & \MM \Z/{\ell^{\nu}}
}
\]
provides a tool for understanding $\MM \Z$ in terms of $\MM \Z/{\ell^{\nu}}$ for each rational prime $\ell$.
Its rational part $\MM \Q$ identifies with its \text{\'et}ale counterpart as in \cite{Voevodsky03b}*{Lemma 6.8}.

\subsection{Algebraic $K$-theory}
\label{subsctn:KGL}

Our next family of examples arises from $K$-theory.
In this section, we will discuss the
motivic analogues of the classical periodic $K$-theory spectrum
$\KU$ and its connective cover $\ku$.

The motivic spectrum $\KGL$ represents algebraic $K$-theory \cite{Voevodsky98} in the following sense.
\index{algebraic $K$-theory spectrum}
\index{$K$-theory spectrum, algebraic}
The homotopy groups $\pi_{n,0} \KGL$ are isomorphic to the algebraic $K$-theory $K_n(k)$ of the base field $k$.
In addition, there is a Bott element $\beta$ in $\pi_{2,1} \KGL$,
and $\KGL$ is periodic in $\beta$.
This determines all of the stable homotopy groups of $\KGL$, in terms of classical algebraic $K$-theory.
The motivic spectrum $\KGL$ is analogous to
$\KU$ in the classical case.

We next wish to consider 
the effective slice cover $\kgl = \f_0 \KGL$, which
is analogous to $\ku$ in the classical case.
In the case of $\KGL$, the very effective slice cover
$\tilde{\f}_0 \KGL$ coincides with $\kgl$ because
$\kgl$ is already very effective \cite{SO}*{Corollary 5.13}.
In the specific $\C$-motivic case, there is another 
way to construct the same motivic spectrum $\kgl$
that is more in the spirit of the classical Postnikov tower \cite{IS11}.

The motivic spectrum $\kgl$ is best understood via its connection
to the motivic truncated Brown-Peterson spectrum
$\MBP \langle 1 \rangle$, discussed below in 
Section \ref{subsctn:MBP}.

\subsection{Motivic Brown-Peterson spectra and truncations}
\label{subsctn:MBP}

Let $\MBP$ denote the motivic Brown-Peterson spectrum at the prime $2$ over a characteristic $0$ base field $k$ constructed in \cite{HuKriz} and \cite{Vezzosi}.
\index{motivic Brown-Peterson spectrum}
\index{Brown-Peterson spectrum, motivic}
This is the universal $2$-typical oriented motivic ring spectrum.
It turns out that the $2$-localized effective slice cover $\kgl_{(2)}$ is the truncated Brown-Peterson spectrum $\MBP\langle{1}\rangle$ sitting in a tower 
\begin{equation*}
\MBP = \MBP\langle{\infty}\rangle
\rightarrow
\cdots
\rightarrow
\MBP\langle{n}\rangle
\rightarrow
\MBP\langle{n-1}\rangle
\rightarrow
\cdots
\rightarrow
\MBP\langle{0}\rangle
\end{equation*}
of $\MBP$-modules.
The maps in this tower come from the fact that $\MBP\langle{n-1}\rangle$
is the cofiber of the map
\[
v_n: \MBP \langle n \rangle \rightarrow
\MBP \langle n \rangle.
\]
This picture is entirely analogous to the classical situation
\cite{JW73},
in which $\BP\langle n - 1 \rangle$ is the cofiber of
\[
\BP \langle n \rangle \rightarrow
\BP \langle n \rangle.
\]

The motivic spectrum
$\MBP\langle{0}\rangle=\MM\Z_{(2)}$ is the $2$-local motivic Eilenberg-Mac Lane spectrum by the theorem of Hopkins, Morel, and Hoyois \cite{hoyois}.
(At odd primes, 
$\MBP\langle{0}\rangle$ is an Adams summand of localized connective algebraic $K$-theory \cite{MR3385689}*{\S4}.)
For $n>1$ we can view the groups $\pi_*\MBP\langle{n}\rangle$ as higher height generalizations of the algebraic $K$-theory groups of the base field.

In order to understand the homology of $\MBP\langle{n}\rangle$ as a comodule over $A_*$ we employ the auxiliary Hopf algebroids
\[\begin{aligned}
  \mathcal{E}(n)_* & =(\M_2,
  A_*/(\xi_1,\xi_2,\ldots)+(\tau_{n+1},\tau_{n+2},\ldots))\\
  &= (\M_2,\M_2[\tau_0,\ldots,\tau_n]/(\tau_i^2-\rho\tau_{i+1},\tau_n^2)).
\end{aligned}\]
We permit $n=\infty$, in which case
\[\begin{aligned}
  \mathcal{E}(\infty)_* & =
  (\M_2,A_*/(\xi_1,\xi_2,\ldots))\\
  &= (\M_2,\M_2[\tau_0,\tau_1,\ldots]/(\tau_i^2 - \rho\tau_{i+1})).
\end{aligned}\]
By  \cite{Ormsby} there is an isomorphism of Hopf algebroids
\[
\MM\Z/2_* \MBP\langle{n}\rangle
\cong
A_*\cotensor_{\mathcal{E}(n)_*} \M_2.
\]
By a standard change-of-rings isomorphism,
the Adams $E_2$-page for $\MBP\langle{n}\rangle$ identifies with 
\[  
\Ext_{\mathcal{E}(n)}(\M_2,\M_2).
\]
For explicit calculations of $\pi_*\MBP\langle{n}\rangle$ over $\C$, $\R$, $\Q_{p}$, and $\Q$ we refer to \cite{MR3073932}.

\subsection{Hermitian $K$-theory}
\label{subsctn:KQ}

In this section, we consider the motivic versions of the
classical real $K$-theory spectrum $\KO$ and its
connective cover $\ko$.

The motivic spectrum $\KQ$ represents Karoubi's hermitian $K$-theory \cite{MR2122220}.
\index{hermitian $K$-theory spectrum}
\index{$K$-theory spectrum, hermitian}
Its effective slice cover $\f_{0}\KQ$ does not coincide with its very effective slice cover $\kq:=\tilde{\f}_{0}\KQ$ studied in \cite{ARO17}.
We note that $\kq$ and $\kgl$ are connected via the motivic Hopf map $\eta\colon\A^{2}\smallsetminus\{0\}\to\PP^{1}$ in the cofiber sequence
\begin{equation}
\label{equation:woodveryeffectivecovers}
\Sigma^{1,1}\kq 
\xrightarrow{\eta} \kq 
\rightarrow 
\kgl.
\end{equation}
An analogous cofiber sequence for the effective slice covers of $\KQ$ and $\KGL$ does not exist by the proof of \cite{rso.hlp}*{Corollary 5.1}.
By \eqref{equation:woodveryeffectivecovers} the Betti realization of $\kq$ identifies with $\ko$ and one can calculate the mod-$2$ motivic homology $\MM\Z/2_{*}\kq$ as 
\[
A_*\cotensor_{A(1)_*} \M_2,
\]
where $A(1)_*$ is the Hopf algebroid
\[\begin{aligned}
  A(1)_* & =(\M_2,
  A_*/(\xi_1^2,\xi_2,\xi_3\ldots)+(\tau_{2},\tau_{3},\ldots))\\
  &= (\M_2,\M_2[\tau_0, \tau_1, \xi_1]/
  (\tau_0^2 + \tau \xi_1 + \rho \tau_1 + \rho \tau_0 \xi_1,
   \xi_1^2, \tau_1^2)
\end{aligned}\]
By change-of-rings,
the Adams $E_2$-page for $\kq$ takes the form
\begin{equation}
\label{equation:masskql}
\Ext^{\ast,*}_{A(1)}(\MM_2,\MM_2).
\end{equation}
As usual, this Adams spectral sequence computes
the homotopy groups of $\kq$ completed at $2$ and $\eta$.
For explicit calculations with \eqref{equation:masskql} 
over $\C$ and $\R$, we refer to \cite{IS11} and \cite{Hill11}.

\subsection{Higher Witt theory}
\label{subsctn:HWT}

In this section we review the proof of Milnor's conjecture on quadratic forms based on the effective slice spectral sequence for higher Witt-theory $\KW$ \cite{roendigs-oestvaer.hermitian}.
Recall that $\KW$ is defined by inverting $\eta$
on $\KQ$.
\index{higher Witt-theory spectrum}
\index{Witt-theory spectrum, higher}

Suppose $k$ is a field of characteristic unequal to $2$.
Recall the Milnor $K$-theory of $k$ 
described in Section \ref{subsctn:Milnor-K}.
In degrees zero, one, and two, these groups agree with Quillen's $K$-groups, but for higher degrees they differ in general. 
Milnor \cite{Milnor69} proposed two conjectures relating $K^{M}_{\ast}(k)/2$ to the mod-$2$ Galois cohomology ring $H^{\ast}(F;\mu_{2})$ and the graded Witt ring 
$GrW^{I}_{\ast}(k)=\oplus_{q\geq 0}I(k)^{q}/I(k)^{q+1}$ for the fundamental ideal $I(k)$ of even dimensional forms
in the Witt ring $W(k)$, 
in the form of two graded ring homomorphisms:
\begin{equation}
\label{equation:milnorconjectures}
\xymatrix{
& K^{M}_{\ast}(k)/2 \ar[dl]_{s^{k}_{\ast}} \ar[dr]^{h^{k}_{\ast}} \\
GrW^{I}_{\ast}(k) && H^{\ast}(k;\mu_{2}).
}
\end{equation}
The Milnor conjecture on Galois cohomology 
states that 
\index{Milnor conjecture}
\index{Galois cohomology}
\index{cohomology, Galois}
$h^{k}_{\ast}$ is an isomorphism, 
while the Milnor conjecture 
on quadratic forms
states that $s^{k}_{\ast}$ is an isomorphism.
The proofs of both conjectures
\cite{Voevodsky03b} \cite{Orlov-Vishik-Voevodsky} 
are two striking applications of motivic homotopy theory.

The slices of $\KW$ were described in 
Part (\ref{part:wittheoryslices}) of Theorem \ref{thm:slice-K}.
We record the first differentials in the effective slice spectral sequence for $\KW$ as maps between motivic spectra.
The differential
\[ 
d_{1}^{\KW}(q)
\colon 
\s_{q}\KW 
\to
\Sigma^{1,0}\s_{q+1}\KW 
\]
is a map of the form
\[ 
\bigvee_{i\in\Z}\Sigma^{2i+q,q}\MM\Z/2 
\to 
\Sigma^{2,1}\bigvee_{j\in\Z}\Sigma^{2j+q,q}\MM\Z/2. 
\]
Let $d_{1}^{\KW}(q,i)$ denote the restriction of $d_{1}^{\KW}(q)$ to the $i$th summand $\Sigma^{2i+q,q}\MM\Z/2$ of $\s_{q}\KW$.
By comparing with motivic cohomology operations of weight one, 
it suffices to consider 
$d_{1}^{\KW}(q,i)$ as a map from $\Sigma^{2i+q,q}\MM\Z/2$ to
\[
\Sigma^{2i+q+4,q+1}\MM\Z/2 
\vee
\Sigma^{2i+q+2,q+1}\MM\Z/2
\vee
\Sigma^{2i+q,q+1}\MM\Z/2. 
\]
The latter map affords the closed formula
\begin{equation}
\label{equation:d1KTdifferentials}
d_{1}^{\KW}(q,i) 
= 
\begin{cases} 
(\Sq^{3}\Sq^{1},\Sq^{2},0) & i-2q\equiv 0\bmod 4 \\
(\Sq^{3}\Sq^{1},\Sq^{2}+\rho\Sq^{1},\tau) & i-2q \equiv 2\bmod 4.
\end{cases} 
\end{equation}

This sets the stage for the proof 
by R{\"o}ndigs and the second author
\cite{roendigs-oestvaer.hermitian}
of Milnor's conjecture on quadratic forms formulated in \cite{Milnor69}*{Question 4.3}.
For fields of characteristic zero this conjecture was first shown by Orlov, Vishik and Voevodsky in \cite{Orlov-Vishik-Voevodsky},
and by Morel \cite{Morel99} using different approaches.

According to Part (\ref{part:wittheoryslices}) 
of Theorem \ref{thm:slice-K}, the effective slice spectral sequence for $\KW$ fills out the entire upper half-plane.
A strenuous computation using \eqref{equation:d1KTdifferentials}, 
Adem relations, 
and the action of the Steenrod squares on the mod-$2$ motivic cohomology ring of $k$ shows that it collapses.
A key point from Theorem \ref{thm:cohlgy-point}
is that if $0 \leq p \leq q$ and 
$a$ belongs to $H^{p,q}(k;\F_2)$, then
$a$ equals $\tau^{q-p} c$ for some
$c$ in $H^{p,p}(k;\F_2)$.
The action of the Steenrod operations in weight at most $1$ is then
completely described by the formulas
\begin{align*}
\Sq^{1}(\tau^nc) 
= &  
\begin{cases} 
\rho\tau^{n-1}c    & n\equiv 1\bmod 2 \\
0                  & n\equiv 0\bmod 2
\end{cases} 
\end{align*}

\begin{align*}
\Sq^{2}(\tau^nc) 
= & 
\begin{cases} 
\rho^{2}\tau^{n-1}c    & n\equiv 2,3\bmod 4\\
0                     & n\equiv 0,1\bmod 4
\end{cases} 
\end{align*}

%
%

With the above in hand, a long calculation shows that
the effective slice $E_2$-page for $\KW$ in degree $(p,q)$ equals
$H^{p,p}(k;\F_2)$ if $p \bmod 4$, and equals zero otherwise.]

To connect this computation with the theory of quadratic forms, 
one shows the spectral sequence converges to the filtration of the Witt ring $W(k)$ by the powers of the fundamental ideal $I(k)$ of even dimensional forms.
By identifying motivic cohomology with Galois cohomology for fields we arrive at the following result.

\begin{theorem} 
\cite{roendigs-oestvaer.hermitian}*{Theorem 1.1}
When $\Char(k)\neq 2$,
the effective slice spectral sequence for $\KW$ converges and furnishes a complete set of invariants 
\[ 
I(k)^{q}/I(k)^{q+1}
\overset{\cong}{\to}
H^{q}(k;\mu_2)
\]
for quadratic forms over $k$ with values in the mod-$2$ Galois cohomology ring.
\end{theorem}

\section{Future directions}
\label{sctn:future}

The purpose of this section is to 
encourage further research into motivic stable
homotopy groups by describing
some specific projects.

\begin{problem}
\index{motivic periodicity}
\index{periodicity, motivic}
\index{motivic nilpotence}
\index{nilpotence, motivic}
Classical periodicity consists of existence results, such as the
existence of Morava $K$-theories and the existence
of $v_n$-self maps; and uniqueness results, such as the
fact that a finite complex possesses a unique $v_n$-self map
(at least up to taking powers).
In the $\C$-motivic situation, many of the analogous
existence results have been established.

On the other hand, much work remains to be done on the uniqueness
results.  In particular, it is not known whether there are any
periodicities in addition to those
discovered by Andrews \cite{Andrews14}, Gheorghe \cite{Gheorghe17}, 
and Krause \cite{Krause18}.
Also, it turns out that a finite complex can possess more than one
type of periodic self-map.  It is not yet understood how periodicities
of different types can co-exist, nor what that means for the structure
of $\C$-motivic stable homotopy groups.

To date, there is little work on exotic nilpotence
and periodicity over other base fields.  At least some of the
same phenomena must occur.  See \cite{Joachimi15} 
and \cite{Hornbostel18} for some first results.
\end{problem}

\begin{problem}
Currently available techniques allow for many more
$\R$-motivic computations than have so far been carried out.
Preliminary data suggests many interesting connections to other
aspects of stable homotopy theory.  In particular, it appears
that the $\R$-motivic image of $J$ has order greater than
the classical image of $J$ \cite{BI19}.
\index{image of $J$}
\index{$J$, image of}
A full accounting of this 
phenomenon is needed.
Both the motivic Adams spectral sequence and the slice
spectral sequence ought to be important tools.

The $\R$-motivic stable homotopy category $\SH(\R)$ appears to be closely
related to the $C_2$-equivariant stable homotopy category.
\index{$C_2$-equivariant stable homotopy}
\index{equivariant stable homotopy, $C_2$}
Data suggests that the $C_2$-equivariant stable homotopy category
is the same as ``$\tau$-periodic $\R$-motivic stable homotopy theory".

Recent on-going work of Behrens and Shah addresses this issue.

Bruner and Greenlees \cite{BG95} showed how 
to reformulate the classical Mahowald invariant \cite{Behrens07}
in terms of $C_2$-equivariant stable homotopy groups.
The $\R$-motivic stable homotopy groups are somewhat easier
to study than the $C_2$-equi\-variant stable homotopy groups,
and they also ought to be useful for Mahowald invariants.
\index{Mahowald invariant}
\index{root invariant}
\end{problem}

\begin{problem}
In classical homotopy, the topological modular forms spectrum
$\tmf$ is an essential tool for studying stable homotopy groups.
\index{topological modular forms}
\index{$\tmf$}
The cohomology of $\tmf$ is $A//A(2)$,
where $A$ is the classical Steenrod algebra and $A(2)$
is the subalgebra of $A$ generated by $\Sq^1$, $\Sq^2$, and $\Sq^4$.
This means that the homotopy groups of $\tmf$
can be computed by an Adams spectral sequence whose 
$E(2)$-page is the cohomology of $A(2)$.
The Adams spectral sequence for $\tmf$ can be completely
analyzed, and this analysis provides much information about
the classical stable homotopy groups.

It is likely that a similar story plays out in the 
motivic situation.  We would like to know that 
a ``motivic modular forms" spectrum $\mmf$ exists over an arbitrary
field $k$.  
\index{motivic modular forms}
\index{$\mmf$}
The motivic stable homotopy groups of $\mmf$ ought to be 
computable by an Adams spectral sequence whose $E_2$-page
is the cohomology of the subalgebra
$A(2)$ of the $k$-motivic Steenrod algebra $A$.
\end{problem}

\begin{problem}
Formula (\ref{equation:sphereslices}) describes the $E_1$-page
of the effective slice spectral sequence, in terms of the 
classical Adams-Novikov $E_2$-page and in terms of
motivic cohomology with coefficients in $\Z$ and $\Z/2^n$.
Preliminary calculations indicate that this
$E_1$-page has interesting ``exotic" products that are not simply
seen by the classical Adams-Novikov $E_2$-page.
This product structure deserves careful study in low dimensions.
See Remark \ref{rmk:slice-mult} for more discussion.
\end{problem}

\begin{problem}
The $E_1$-page of the effective slice spectral sequence can be 
described over a general base field $k$, at least in terms of
the Milnor $K$-theory of $k$.  Some differentials have been
understood in general in very low dimensions
\cite{roendigs-spitzweck-oestvaer.1line}.
These results about effective slice differentials ought to be
accessible in a larger range of dimensions, at least relative
to arithmetic input from $k$.
\end{problem}

\begin{problem}
Section \ref{subsctn:HWT} discusses how the effective slice filtration
for $\KW$ informs us about quadratic forms over a field $k$.
It is possible that these ideas can be extended to study
quadratic forms over smooth $k$-schemes.
\end{problem}

\bibliographystyle{amsalpha}

\bibliographystyle{amsalpha}

\end{document}